# GLOBAL DYNAMICAL SOLVERS FOR NONLINEAR PROGRAMMING PROBLEMS

**Iasson Karafyllis[*] and Miroslav Krstic[**]**

[*]Dept. of Mathematics, National Technical University of Athens, 15780, Athens, Greece, email: iasonkar@central.ntua.gr

[**]Dept. of Mechanical and Aerospace Eng., University of California, San Diego, La Jolla, CA 92093-0411, U.S.A., email: krstic@ucsd.edu

**Abstract**

We construct a family of globally defined dynamical systems for a nonlinear programming problem, such that: (a) the equilibrium points are the unknown (and sought) critical points of the problem, (b) for every initial condition, the solution of the corresponding initial value problem converges to the set of critical points, (c) every strict local minimum is locally asymptotically stable, (d) the feasible set is a positively invariant set, and (e) the dynamical system is given explicitly and does not involve the unknown critical points of the problem. No convexity assumption is employed. The construction of the family of dynamical systems is based on an extension of the Control Lyapunov Function methodology, which employs extensions of LaSalle's theorem and are of independent interest. Examples illustrate the obtained results.

## 1. Introduction

Dynamical systems have been used in the past for the solution of Nonlinear Programming (NLP) problems. The reader may consult [2,8,9,13,24,29,34,35] for various results on the topic. Some methods are interior-point methods (in the sense that are defined only on the feasible set) while other methods are exterior-point methods (in the sense that are defined at least in a neighborhood of the feasible set). As remarked in [7,14,18], each system of ordinary differential equations that solves a NLP problem when combined with a numerical scheme for solving Ordinary Differential Equations (ODEs) provides a numerical scheme for solving the NLP problem. Dynamical systems have also been utilized for the solution of Linear Programming and NLP problems in the literature of neural networks (see for example [31,32,33] as well as the review paper [30] and the references therein).

Therefore, it is justified to use the term *"dynamical NLP solver"* for a dynamical system for which some of its solutions converge to the solutions of a NLP problem. The recent work [20] applied feedback stabilization methods for the explicit construction of interior-point dynamical NLP solvers. However, interior-point dynamical NLP solvers have some disadvantages:

(a) they have to be initiated in the feasible set, and

(b) the application of a numerical integrator is problematic since the system is defined only on the feasible set $S \subseteq \Re^n$. Thus, the numerical integration may involve a projection on the feasible set (a serious complication; see [15,16]).

In this work, we are interested in the application of feedback stabilization methods for constructing *exterior-point* dynamical NLP solvers. More specifically, we consider a standard NLP problem with sufficient regularity properties and so that necessary Karush-Kuhn-Tucker (KKT) conditions of the NLP hold. Inspired by the methods employed in the book [17], our goal is to construct a globally defined dynamical system with the following properties:

Property 1: The vector field appearing on the right hand side of the dynamical system is a locally Lipschitz vector field which is globally defined. This property is required for uniqueness of the solutions of the dynamical system. Moreover, this property is required because we would like to be able to apply Runge-Kutta schemes for the simulation of the solutions of the dynamical system.

Property 2: The equilibrium points of the dynamical system are exactly the points for which the necessary Karush-Kuhn-Tucker conditions of the NLP hold.

Property 3: The vector field appearing on the right hand side of the dynamical system must be explicitly known. Formulas for the vector field must be provided: the formulas must not involve the solution of the NLP problem.

Property 4: For every initial point, the solution of the corresponding initial value problem converges to the set of KKT points. Moreover, every strict local minimum which is an isolated KKT point is locally asymptotically stable.

Property 5: The feasible set is a positively invariant set for the dynamical system. This property may be important for certain applications.

Property 6: All previous properties must be valid for general NLP problems without any convexity assumption.

It must be noted that the properties 1-6 are *rarely satisfied* by other differential equation methods for solving NLPs. For example, in [2] and [8], dynamical NLP solvers are proposed for certain NLP problems. However, the solution of the NLP problem is not an equilibrium point for the constructed time-varying dynamical system in [8]. Antipin in [2] constructs an autonomous dynamical system for which the solution of the NLP problem is an equilibrium point and for which the locally Lipschitz vector field appearing in the right hand side of the dynamical system does not depend on the location of the unknown point. However, the definition of the vector field appearing on the right hand side of the dynamical system is involved (it requires the solution of a NLP since it involves a projection on the feasible set). Special NLP problems under additional convexity assumptions have been studied in [31]. Convexity assumptions appear in almost all neural networks proposed for the solution of mathematical programming problems (see [32,33] and the references in the review paper [30]). On the other hand, the papers [29,34] propose systems of differential equations that satisfy properties 1-6 for systems without inequality constraints. Local results are provided in the paper [35] and differential equations based on barrier methods were considered in [9].

The feedback stabilization method employed in this work is the Control Lyapunov Function methodology (CLF; see [3,11,19,27]). However, we face the important issue of the construction of a CLF which must combine in an appropriate way a penalty term $V(x)$ (i.e., a function which is zero on the feasible set and positive out of the feasible set) and the objective function $\theta(x)$ (which is a natural candidate for the Lyapunov function on the feasible set). Such a combination is very difficult and can be achieved under very demanding assumptions.

In order to overcome the Lyapunov construction problem we propose the idea of using two functions as Lyapunov-like functions: the penalty term $V(x)$ when we are away from the feasible set and the objective function $\theta(x)$ when we are on the feasible set. Moreover, we don't use a feedback construction methodology which is based on the Lyapunov theorem. Instead, the feedback construction in this work is based on our extensions of the LaSalle's theorem, which are of independent interest. Therefore, the contribution of the paper is threefold:

1) Extensions of LaSalle's theorem are provided.
2) The solution of a special feedback stabilization problem is presented. The provided solution is based on an extension of the CLF methodology, which employs the obtained extensions of LaSalle's theorem.
3) Dynamical NLP solvers with the aforementioned properties 1-6 are constructed, based on the solution of the special feedback stabilization problem mentioned above.

The construction of the dynamical NLP solvers with the aforementioned properties 1-6 involve the linear independence constraint qualification. The linear independence constraint qualification which is assumed in this work is a restrictive assumption: it is more restrictive than the Mangasarian-Fromovitz constraint qualification in [23] or the constant rank constraint qualification (see [1] and references therein), which are all more restrictive than the Guignard constraint qualification (see [4]). However, the linear independence constraint qualification has the advantage of being easily checkable and of being true in many interesting cases (the work [25] showed that this assumption holds generically) and it is a vital ingredient for many numerical methods (successive quadratic programming-see [10,26]). Furthermore, the linear independence constraint qualification allows us to obtain easy formulas for the required vector field.

The structure of the paper is as follows: Section 2 describes the problem studied in this paper. Section 3 contains the extensions of LaSalle's theorem, while Section 4 provides the solution of certain feedback stabilization problems based on the obtained extensions of LaSalle's theorem. The feedback stabilization problem studied in Section 4 is a special problem, which can be used for the construction of dynamical NLP solvers. Section 5 is devoted to the construction of the dynamical NLP solvers, based on the results of the previous section. Special cases for NLP problems for which the formulas of the dynamical NLP solver become simpler are presented in Section 6. Section 7 contains three illustrative examples and Section 8 provides the concluding remarks of the present work. The Appendix contains the proofs of certain auxiliary results.

*Notation.* Throughout the paper we adopt the following notation:

* Let $A \subseteq \Re^n$ be an open set. By $C^0(A;\Omega)$, we denote the class of continuous functions on $A$, which take values in $\Omega$. By $C^k(A;\Omega)$, where $k \geq 1$ is an integer, we denote the class of differentiable functions on $A$ with continuous derivatives up to order $k$, which take values in $\Omega$. By $C^\infty(A;\Omega)$, we denote the class of differentiable functions on $A$ having continuous derivatives of all orders (smooth functions), which take values in $\Omega$, i.e., $C^\infty(A;\Omega) = \bigcap_{k \geq 1} C^k(A;\Omega)$.

* For a vector $x \in \Re^n$ we denote by $|x|$ its usual Euclidean norm and by $x'$ its transpose. For a real matrix $A \in \Re^{n \times m}$ we denote by $A' \in \Re^{m \times n}$ its transpose. $I_n \in \Re^{n \times n}$ denotes the identity matrix. For a square matrix $A \in \Re^{n \times n}$, $\det(A)$ denotes its determinant. The adjugate matrix $adj(A)$ of a square matrix $A \in \Re^{n \times n}$ is the transpose of its cofactor matrix, i.e., it is the square matrix that satisfies $A\,adj(A) = adj(A)A = \det(A)I_n$. For $x \in \Re^n$, $diag(x) \in \Re^{n \times n}$ denotes the diagonal matrix with $x \in \Re^n$ on its diagonal.

* For every $x = (x_1,...,x_n)' \in \Re^n$ we define $x^+ = \big(\max(0,x_1),...,\max(0,x_n)\big)' \in \Re^n$ and $x^- = \big(\min(0,x_1),...,\min(0,x_n)\big)' \in \Re^n$. Notice that the following property holds for every positive definite and diagonal matrix $R \in \Re^{n \times n}$: $x'Rx^+ = 0 \Leftrightarrow x^+ = 0$.

* $\Re_+^n := (\Re_+)^n = \{(x_1,...,x_n)' \in \Re^n : x_1 \geq 0,...,x_n \geq 0\}$. Let $x, y \in \Re^n$. We say that $x \leq y$ if and only if $(y-x) \in \Re_+^n$.

* $K_\infty$ is the class of continuous, increasing, unbounded functions $a: \Re_+ \to \Re_+$ with $a(0) = 0$.

* For every scalar continuously differentiable function $V: \Re^n \to \Re$, $\nabla V(x)$ denotes the gradient of $V$ at $x \in \Re^n$, i.e., $\nabla V(x) = \left(\frac{\partial V}{\partial x_1}(x),...,\frac{\partial V}{\partial x_n}(x)\right)$.

## 2. Problem Description

Consider the Nonlinear Programming problem:
$$\min\{\theta(x) : x \in S\} \quad (2.1)$$

where $x \in \Re^n$ and the closed set $S \subseteq \Re^n$ is defined by
$$S := \left\{ x \in \Re^n : h_1(x) = ... = h_m(x) = 0, \max_{j=1,...k}\big(g_j(x)\big) \leq 0 \right\} \quad (2.2)$$

where $m < n$. All mappings $\theta: \Re^n \to \Re$, $h_i: \Re^n \to \Re$ ($i = 1,...,m$), $g_j: \Re^n \to \Re$ ($j = 1,...,k$) are twice continuously differentiable.

**(H1)** *The feasible set $S \subseteq \Re^n$ defined by (2.2) is non-empty and the sublevel sets of $\theta: \Re^n \to \Re$ are compact sets, i.e., for every $x_0 \in S$ the level set*
$$\{x \in S : \theta(x) \leq \theta(x_0)\}$$
*is compact.*

Assumption (H1) is a standard assumption which guarantees that the NLP problem described by (2.1) and (2.2) is well-posed and admits at least one global solution, $x^* \in S$.

We define:
$$h(x) = \begin{bmatrix} h_1(x) \\ \vdots \\ h_m(x) \end{bmatrix} \in \Re^m, \ A(x) = \begin{bmatrix} \nabla h_1(x) \\ \vdots \\ \nabla h_m(x) \end{bmatrix} \in \Re^{m \times n},$$

$$g(x) = \begin{bmatrix} g_1(x) \\ \vdots \\ g_k(x) \end{bmatrix} \in \Re^k, \ B(x) = \begin{bmatrix} \nabla g_1(x) \\ \vdots \\ \nabla g_k(x) \end{bmatrix} \in \Re^{k \times n}, \text{ for all } x \in \Re^n \quad (2.3)$$

We also assume the "Mangasarian-Fromovitz constraint qualification".

**(H2)** *The row vectors $\nabla h_i(x)$ ($i=1,...,m$) are linearly independent for all $x \in S$. Moreover, for every $x \in S$, if $\{i \in \{1,...,k\}: g_i(x)=0\} \neq \emptyset$ then there exists $\xi \in \Re^n$ such that $\nabla g_j(x)\xi < 0$ for all $j=1,...,k$ for which $g_j(x)=0$ (active constraints) and $\nabla h_i(x)\xi = 0$ for all $i=1,...,m$.*

Next, we define the set of critical points for the NLP problem defined by (2.1) and (2.2).

Let $\Phi \subseteq S$ be the set of all points $x \in S$ for which there exist $\lambda \in \Re^m$ and $\mu \in \Re^k_+$ such that the following equations hold:

$$\begin{aligned}(\nabla \theta(x))' + A'(x)\lambda + B'(x)\mu &= 0 \\ \mu' g(x) &= 0\end{aligned} \quad (2.4)$$

In other words, $\Phi \subseteq S$ is the set of critical points or Karush-Kuhn-Tucker (KKT) points for the problem defined by (2.1) and (2.2). The "Mangasarian-Fromovitz constraint qualification" (Assumption (H2)) guarantees that every solution of (2.1) and (2.2) is a KKT point.

As described in the Introduction, the problem studied in the present paper is the construction of a globally defined locally Lipschitz vector field $f(x)$ such that the solutions of the dynamical system $\dot{x} = f(x)$, $x \in \Re^n$ converge to $\Phi$ for all initial conditions. The construction of the vector field should not involve the set $\Phi \subseteq S$ itself, because the set $\Phi \subseteq S$ is unknown (it is what we want to determine).

## 3. Extensions of LaSalle's Theorem

While LaSalle's theorem deals with one function $V \in C^1(\Re^n; \Re)$ so that the set $\{x \in \Re^n : \nabla V(x)f(x) = 0\}$ is a global attractor for the dynamical system $\dot{x} = f(x)$, $x \in \Re^n$, in this section we present conditions for two functions $V \in C^1(\Re^n; \Re)$, $\theta \in C^1(\Re^n; \Re)$ for which the set $\{x \in \Re^n : \nabla \theta(x)f(x) = \nabla V(x)f(x) = 0\}$ is a global attractor. We start by presenting conditions for the case of weak attractor. The proof of Theorem 3.1 is provided in the Appendix.

**Theorem 3.1 (First Extension of LaSalle's theorem-the weak attractor case):** *Let $f : \Re^n \to \Re^n$ be a locally Lipschitz vector field and let $V \in C^1(\Re^n; \Re)$, $\theta \in C^1(\Re^n; \Re)$ be functions that satisfy:*

$$\nabla V(x)f(x) \leq 0, \text{ for all } x \in \Re^n \quad (3.1)$$

$$\nabla \theta(x)f(x) \leq 0, \text{ for all } x \in \Re^n \text{ with } \nabla V(x)f(x) = 0 \quad (3.2)$$

*Suppose that for every $z \in \Re^n$ and for every $y \in \{x \in \Re^n : V(x) \leq V(z)\}$ the set $\{x \in \Re^n : V(x) \leq V(z), \theta(x) \leq \theta(y)\}$ is compact. Moreover, assume that for every $z \in \Re^n$ there exists $y_z \in \{x \in \Re^n : V(x) \leq V(z)\}$ with $\theta(y_z) \geq \theta(z)$ such that*

$$\{x \in \Re^n : V(x) \leq V(z)\} \subseteq \{x \in \Re^n : \theta(x) \leq \theta(y_z)\} \cup \{x \in \Re^n : \nabla \theta(x)f(x) \leq 0\} \quad (3.3)$$

*Consider the dynamical system*

$$\dot{x} = f(x), \, x \in \Re^n \tag{3.4}$$

*Then for every $x_0 \in \Re^n$ the unique solution $x(t)$ of the initial value problem (3.4) with $x(0) = x_0$ is defined for all $t \geq 0$ and is bounded. If we further denote by $\omega(x_0)$ the positive limit set of $\{x(t): t \geq 0\}$, then for every $x_0 \in \Re^n$, $\omega(x_0)$ is a non-empty, compact, connected, invariant set which satisfies $\omega(x_0) \subseteq \{x \in \Re^n : \nabla V(x) f(x) = 0\}$ and $\omega(x_0) \cap \{x \in \Re^n : \nabla \theta(x) f(x) = \nabla V(x) f(x) = 0\} \neq \varnothing$.*

**Remark 3.2:** (a) Using the terminology in [6], if the set $M = \{x \in \Re^n : \nabla \theta(x) f(x) = \nabla V(x) f(x) = 0\}$ is compact then it is a weak attractor for the dynamical system (3.4) with region of weak attraction being the whole space $\Re^n$.
(b) If the set $\{x \in \Re^n : V(x) \leq V(z)\}$ is bounded then there exists $y_z \in \{x \in \Re^n : V(x) \leq V(z)\}$ with $\theta(y_z) \geq \theta(z)$ such that (3.3) holds (select $y_z \in \Re^n$ to satisfy $\theta(y_z) = \max\{\theta(x) : x \in \Re^n, V(x) \leq V(z)\}$). However, inclusion (3.3) may hold even for cases where the set $\{x \in \Re^n : V(x) \leq V(z)\}$ is not bounded.

Next we present conditions for the case of global attractor. The proof of Theorem 3.3 is provided in the Appendix.

**Theorem 3.3 (Second Extension of LaSalle's theorem-the global attractor case):** *Suppose that the assumptions of Theorem 3.1 hold. Furthermore, suppose that $V(x) \geq 0$ for all $x \in \Re^n$ and that the following equation holds:*

$$S = \{x \in \Re^n : \nabla V(x) f(x) = 0\} = \{x \in \Re^n : V(x) = 0\} \tag{3.5}$$

*Finally, suppose that for every compact set $K \subseteq \Re^n$ with $K \cap S \neq \varnothing$ there exist a function $g \in K_\infty \cap C^1((0, +\infty); \Re_+)$ and a constant $\delta > 0$ with the following property*

$$\nabla \theta(x) f(x) \leq \frac{dg}{ds}(V(x)) |\nabla V(x) f(x)|, \text{ for all } x \in K \text{ with } 0 < V(x) \leq \delta \text{ and } |\nabla V(x) f(x)| \leq \delta \tag{3.6}$$

*Then for every $x_0 \in \Re^n$ the unique solution $x(t)$ of the initial value problem (3.4) with $x(0) = x_0$ is defined for all $t \geq 0$ and is bounded. The set $S$, defined in (3.5), is a positively invariant set for the dynamical system (3.4). If we further denote by $\omega(x_0)$ the positive limit set of $\{x(t): t \geq 0\}$, then for every $x_0 \in \Re^n$, $\omega(x_0)$ is a non-empty, compact, connected, invariant set which satisfies $\omega(x_0) \subseteq \{x \in \Re^n : V(x) = \nabla \theta(x) f(x) = 0\}$. Moreover, every equilibrium point $x^* \in S$ of the dynamical system (3.4), which satisfies $\theta(x^*) < \theta(x)$ for all $x \in S \setminus \{x^*\}$ with $|x - x^*| < \tilde{\delta}$ and $\{x \in S : \nabla \theta(x) f(x) = 0, |x - x^*| < \tilde{\delta}\} = \{x^*\}$, for an appropriate constant $\tilde{\delta} > 0$, is a locally asymptotically stable equilibrium point of the dynamical system (3.4).*

The following example illustrates the use of Theorem 3.3 for the analysis of the qualitative behavior of a dynamical system.

**Example 3.4:** Consider the nonlinear planar system (3.4) with

$$f(x) = -Q(x)\left[\left(\max\left(0, x_1^2 + x_2^2 - 1\right)\right)^2 + \frac{Q(x)}{4}\right]\begin{bmatrix}1\\1\end{bmatrix}$$
$$+ \left(4\left(\max\left(0, x_1^2 + x_2^2 - 1\right)\right)^2(x_1 + x_2 - 1) + (x_1 + x_2)Q(x) - \max(0, x_1 + x_2)\right)\begin{bmatrix}x_1\\x_2\end{bmatrix}, \text{ for all } x \in \Re^2 \quad (3.7)$$

where

$$Q(x) := 4(x_1^2 + x_2^2) - \min\left(0, x_1^2 + x_2^2 - 1\right), \quad (3.8)$$

Next, we use Theorem 3.3 with

$$V(x) = \frac{1}{2}\left(\max\left(0, x_1^2 + x_2^2 - 1\right)\right)^2, \quad \theta(x) = x_1 + x_2 \quad (3.9)$$

Straightforward calculations reveal that the following equations hold for all $x \in \Re^2$:

$$\nabla V(x) f(x) = -2\max\left(0, x_1^2 + x_2^2 - 1\right)\left(x_1^2 + x_2^2\right)\left(4\left(\max\left(0, x_1^2 + x_2^2 - 1\right)\right)^2 + \max(0, x_1 + x_2)\right) \quad (3.10)$$

$$\nabla \theta(x) f(x) = -4(x_1 - x_2)^2\left(\max\left(0, x_1^2 + x_2^2 - 1\right)\right)^2 - (x_1 + x_2)\max(0, x_1 + x_2)$$
$$-\left((x_1 - x_2)^2 - \frac{\min\left(0, x_1^2 + x_2^2 - 1\right)}{2}\right)Q(x) - 4\left(\max\left(0, x_1^2 + x_2^2 - 1\right)\right)^2(x_1 + x_2) \quad (3.11)$$

Equation (3.10) shows that inequality (3.1) and equation (3.5) hold. Moreover, definition (3.9) shows that for every $z \in \Re^2$ the set $\{x \in \Re^2 : V(x) \leq V(z)\}$ is compact. Therefore, Remark 3.2(b) implies that for every $z \in \Re^2$ there exists $y_z \in \{x \in \Re^2 : V(x) \leq V(z)\}$ with $\theta(y_z) \geq \theta(z)$ such that inclusion (3.3) holds. Furthermore, equation (3.11) and definition (3.9) shows that when $V(x) = 0$, we have:

$$\nabla \theta(x) f(x) = -\left((x_1 - x_2)^2 - \frac{\min\left(0, x_1^2 + x_2^2 - 1\right)}{2}\right)Q(x) - (x_1 + x_2)\max(0, x_1 + x_2) \quad (3.12)$$

Equation (3.12) shows that inequality (3.2) holds. Finally, we show that for every compact set $K \subseteq \Re^2$ with $K \cap S \neq \emptyset$ there exist a function $g \in K_\infty \cap C^1((0, +\infty); \Re_+)$ and a constant $\delta > 0$ satisfying (3.6). Let $K \subseteq \Re^2$ be a compact set with $K \cap S \neq \emptyset$. Equation (3.11) implies that

$$\nabla \theta(x) f(x) \leq 4M\left(\max\left(0, x_1^2 + x_2^2 - 1\right)\right)^2, \text{ for all } x \in K \quad (3.13)$$

with $M = \max\{x_1 + x_2 : x \in K\}$. Equation (3.10) and definition (3.9) show that

$$|\nabla V(x) f(x)| \geq 8\left(\max\left(0, x_1^2 + x_2^2 - 1\right)\right)^3 \text{ for all } x \in \Re^2 \text{ with } V(x) > 0. \quad (3.14)$$

Inequalities (3.13), (3.14) and definition (3.9) show that the inequality $\nabla \theta(x) f(x) \le \frac{M}{2\sqrt{2V(x)}} |\nabla V(x) f(x)|$ holds for all $x \in K$ with $V(x) > 0$. Consequently, (3.6) holds with arbitrary $\delta > 0$ and $g(s) := \frac{M}{\sqrt{2}} \sqrt{s}$, which is a function of class $g \in K_\infty \cap C^1((0,+\infty); \Re_+)$.

We conclude that all assumptions of Theorem 3.3 hold. Definition (3.9) and equation (3.12) imply that $\{x \in \Re^2 : V(x) = \nabla \theta(x) f(x) = 0\} = \left\{\left(-\frac{\sqrt{2}}{2}, -\frac{\sqrt{2}}{2}\right)'\right\}$ and consequently, Theorem 3.3 implies that for every $x_0 \in \Re^n$ the unique solution $x(t)$ of the initial value problem (3.4) with $x(0) = x_0$ is defined for all $t \ge 0$ and satisfies $\lim_{t \to +\infty} x(t) = x^* = \left(-\frac{\sqrt{2}}{2}, -\frac{\sqrt{2}}{2}\right)'$ (global attractivity). Finally, since $\theta\left(-\frac{\sqrt{2}}{2}, -\frac{\sqrt{2}}{2}\right) < \theta(x)$ for all $x \in S \setminus \{x^*\}$, Theorem 3.3 allows us to conclude that $\left(-\frac{\sqrt{2}}{2}, -\frac{\sqrt{2}}{2}\right)'$ is locally asymptotically stable, which combined with global attractivity, implies that $\left(-\frac{\sqrt{2}}{2}, -\frac{\sqrt{2}}{2}\right)'$ is globally asymptotically stable.

The obtained result would be difficult to obtain by a single Lyapunov-like function, since we cannot easily construct a Lyapunov function for the dynamical system (3.4) with (3.7), (3.8). ◁

## 4. Feedback Construction

In this section we provide solutions to the following problem. We are given two functions $V \in C^1(\Re^n; \Re_+)$, $\theta \in C^1(\Re^n; \Re)$ for which the following equality holds

$$S = \left\{ x \in \Re^n : |\nabla V(x)| = 0 \right\} = \left\{ x \in \Re^n : V(x) = 0 \right\} \quad (4.1)$$

We are also given a vector field $F(x)$ defined on a neighborhood of the set $S$, which satisfies the inequalities $\nabla V(x) F(x) \le 0$, $\nabla \theta(x) F(x) \le 0$ for all $x \in \Re^n$ in a neighborhood of $S$. Our goal is the explicit design of a locally Lipschitz feedback law $u = f(x)$ for the control system $\dot{x} = u \in \Re^n$ so that the set $\{x \in S : \nabla \theta(x) F(x) = 0\}$ is a global attractor for the closed-loop system. This is a special feedback stabilization problem and the reader may wonder why such a problem is studied. However, the following section shows that this special feedback stabilization problem is exactly the problem needed to be solved when dealing with the construction of a dynamical NLP solver.

The following theorem provides a solution of the above feedback control problem, which provides an explicit formula for the locally Lipschitz feedback law $u = f(x)$ and is based on Theorem 3.3.

**Theorem 4.1:** *Let $V \in C^1(\Re^n; \Re_+)$, $\theta \in C^1(\Re^n; \Re)$ be given functions with locally Lipschitz partial derivatives and suppose that (4.1) holds. Suppose that for every $z \in \Re^n$ and for every $y \in \{x \in \Re^n : V(x) \le V(z)\}$ the set $\{x \in \Re^n : V(x) \le V(z), \theta(x) \le \theta(y)\}$ is compact. Let $\Omega \in C^0(\Re^n; \Re_+)$ be a* *locally Lipschitz function with* $\Omega(x) > 0$ *for* $x \notin S$, $\Omega(x) = 0$ *for all* $x \in S$. *Suppose that there exist a function* $\gamma \in C^0(\tilde{S}; (0, +\infty))$, *where* $\tilde{S} = \{x \in \Re^n : \Omega(x) < 1\}$, *and a locally Lipschitz vector field* $F : \tilde{S} \to \Re^n$ *such that*

$$\gamma(x)|\nabla V(x)|^2 \geq V(x) \text{ for all } x \in \Re^n \text{ with } \Omega(x) < 1 \quad (4.2)$$

$$\nabla V(x) F(x) \leq 0, \ \nabla \theta(x) F(x) \leq 0, \text{ for all } x \in \Re^n \text{ with } \Omega(x) < 1 \quad (4.3)$$

*Finally, suppose that there exists a locally Lipschitz function* $\psi : \Re^n \to (0, +\infty)$ *for which the following property holds:*

(*) *For every* $z \in \Re^n$ *there exists* $y_z \in \{x \in \Re^n : V(x) \leq V(z)\}$ *with* $\theta(y_z) \geq \theta(z)$ *such that*

$$\{x \in \Re^n : a(x) > 0, V(x) \leq V(z)\} \subseteq \{x \in \Re^n : \theta(x) \leq \theta(y_z)\} \quad (4.4)$$

*where*

$$a(x) := -\psi(x)\nabla V(x)(\nabla \theta(x))' - |\nabla V(x)|^2 |\nabla \theta(x)|^2 + |\nabla V(x)(\nabla \theta(x))'|^2 \quad (4.5)$$

*Let* $\beta > 1$ *be a constant,* $\sigma : \Re^n \to (0, +\infty)$ *be a locally Lipschitz function and define the locally Lipschitz vector field:*

$$f(x) := \sigma(x)(1 - \beta\Omega(x))F(x) - \beta\sigma(x)\Omega(x)\psi(x)(\nabla V(x))' - \beta\sigma(x)\Omega(x)\left(|\nabla V(x)|^2 I_n - (\nabla V(x))'\nabla V(x)\right)(\nabla \theta(x))',$$

*for all* $x \in \Re^n$ *with* $\beta\Omega(x) < 1$ \quad (4.6)

$$f(x) := -\sigma(x)\psi(x)(\nabla V(x))' - \sigma(x)\left(|\nabla V(x)|^2 I_n - (\nabla V(x))'\nabla V(x)\right)(\nabla \theta(x))', \text{ for all } x \in \Re^n \text{ with } \beta\Omega(x) \geq 1 \quad (4.7)$$

*Then for every* $x_0 \in \Re^n$ *the unique solution* $x(t)$ *of the initial value problem (3.4) with* $x(0) = x_0$ *is defined for all* $t \geq 0$ *and is bounded. The set* $S$, *defined in (4.1), is a positively invariant set for the dynamical system (3.4). If we further denote by* $\omega(x_0)$ *the positive limit set of* $\{x(t) : t \geq 0\}$, *then for every* $x_0 \in \Re^n$, $\omega(x_0)$ *is a non-empty, compact, connected, invariant set which satisfies* $\omega(x_0) \subseteq \{x \in S : \nabla \theta(x)F(x) = 0\}$. *Moreover, every point* $x^* \in S$, *which satisfies* $F(x^*) = 0$, $\theta(x^*) < \theta(x)$ *for all* $x \in S \setminus \{x^*\}$ *with* $|x - x^*| < \tilde{\delta}$ *and* $\{x \in S : \nabla\theta(x)F(x) = 0, |x - x^*| < \tilde{\delta}\} = \{x^*\}$, *for an appropriate constant* $\tilde{\delta} > 0$, *is a locally asymptotically stable equilibrium point of the dynamical system (3.4).*

**Proof:** It suffices to show that all assumptions of Theorem 3.3 hold for the vector field $f : \Re^n \to \Re^n$ defined by (4.6) and (4.7). Indeed, definitions (4.5), (4.6), (4.7) imply that

$$\nabla V(x) f(x) = \sigma(x)(1 - \beta\Omega(x))\nabla V(x)F(x) - \beta\sigma(x)\Omega(x)\psi(x)|\nabla V(x)|^2 \leq 0,$$

for all $x \in \Re^n$ with $\beta\Omega(x) < 1$ \quad (4.8)

$$\nabla V(x)f(x) = -\sigma(x)\psi(x)|\nabla V(x)|^2 \leq 0, \text{ for all } x \in \Re^n \text{ with } \beta\Omega(x) \geq 1 \quad (4.9)$$

$$\nabla\theta(x)f(x) = \sigma(x)(1 - \beta\Omega(x))\nabla\theta(x)F(x) + \beta\sigma(x)\Omega(x)a(x),$$

for all $x \in \Re^n$ with $\beta\Omega(x) < 1$ \quad (4.10)

$$\nabla\theta(x)f(x) = \sigma(x)a(x), \text{ for all } x \in \Re^n \text{ with } \beta\Omega(x) \geq 1 \quad (4.11)$$

More specifically, inequality (4.8) is a consequence of (4.3). Consequently, inequality (3.1) holds and it holds that $\nabla V(x)f(x)=0$ only when $x\in S$. The latter fact and definition (4.6) (which shows that $f(x)=\sigma(x)F(x)$ when $x\in S$) in conjunction with (4.3), imply that inequality (3.2) holds. It follows from (4.3), (4.10), (4.11) that the following implication holds:

$$\text{If } a(x)\leq 0 \text{ then } \nabla\theta(x)f(x)\leq 0 \tag{4.12}$$

Property (*) and implication (4.12) guarantee the inclusion (3.3).

Finally, we show that for every compact set $K\subseteq\Re^n$ with $K\cap S\neq\emptyset$ there exist constants $M,\delta>0$ with the following property

$$\nabla\theta(x)f(x)\leq\frac{M}{\sqrt{V(x)}}|\nabla V(x)f(x)|, \text{ for all } x\in K \text{ with } 0<V(x)\leq\delta \text{ and } |\nabla V(x)f(x)|\leq\delta \tag{4.13}$$

In other words, we show that inequality (3.6) holds with $g(s):=2M\sqrt{s}$ for $s\geq 0$. In order to show the validity of (4.13), we need the following claim.

**Claim:** *There exist functions* $\tilde{\gamma}\in C^0(\Re^n;(0,+\infty))$, $p\in K_\infty$ *such that*

$$\tilde{\gamma}(x)p(|\nabla V(x)|)\geq\Omega(x) \text{ for all } x\in\Re^n \tag{4.14}$$

**Proof of Claim:** Define the function $q:\Re_+\to[0,1]$ by means of the formula

$$q(s):=\sup\left\{\frac{\Omega(x)}{(1+\Omega(x))(1+|x|^2)}:|\nabla V(x)|\leq s\right\}, \text{ for all } s\geq 0 \tag{4.15}$$

Since the function $\frac{\Omega(x)}{(1+\Omega(x))(1+|x|^2)}$ is non-negative and bounded by 1, it follows that $q(s)$ is well-defined by (4.15) and satisfies $q(s)\in[0,1]$ for all $s\geq 0$. Notice that since $\Omega(x)=0$ for all $x\in S$ and since (4.1) holds, definition (4.15) implies that $q(0)=0$. Moreover, $q:\Re_+\to[0,1]$ is a non-decreasing function, which satisfies the following inequality for all $x\in\Re^n$:

$$(1+|x|^2)(1+\Omega(x))q(|\nabla V(x)|)\geq\Omega(x), \text{ for all } x\in\Re^n \tag{4.16}$$

Next, we show that $\lim_{s\to 0^+}(q(s))=q(0)=0$. It suffices to show that $\limsup_{s\to 0^+}(q(s))=0$. Suppose, on the contrary that $\limsup_{s\to 0^+}(q(s))=l>0$. Then there exists a sequence $\{s_i>0\}_{i=0}^\infty$ with $s_i\to 0$ and $q(s_i)\geq l/2$. Consequently, definition (4.15) implies that there exists a sequence $\{x_i\in\Re^n\}_{i=0}^\infty$ with $|\nabla V(x_i)|\to 0$ and $\frac{\Omega(x_i)}{(1+\Omega(x_i))(1+|x_i|^2)}\geq l/4$. The inequality $\frac{\Omega(x_i)}{(1+\Omega(x_i))(1+|x_i|^2)}\geq l/4$ implies the inequality $4l^{-1}-1\geq|x_i|^2$, which shows that the sequence $\{x_i\in\Re^n\}_{i=0}^\infty$ is bounded. Consequently, there exists a subsequence still denoted by $\{x_i\in\Re^n\}_{i=0}^\infty$, which converges, i.e., there exists $x^*\in\Re^n$ with $x_i\to x^*$. By continuity, we have $|\nabla V(x^*)|=0$ and $\frac{\Omega(x^*)}{(1+\Omega(x^*))(1+|x^*|^2)}\geq l/4$. Since $\Omega(x)=0$ for all $x\in S$ and since (4.1) holds we get $\Omega(x^*)=0$, which combined with the inequality $\dfrac{\Omega(x^*)}{(1+\Omega(x^*))(1+|x^*|^2)} \geq l/4$, contradicts the assumption that $\limsup_{s\to 0^+}(q(s))=l>0$. Therefore, $\limsup_{s\to 0^+}(q(s))=0$.

Lemma 2.4 on page 65 in [19] implies that there exists $p \in K_\infty$ such that $q(s) \leq p(s)$ for all $s \geq 0$. Inequality (4.14) is a direct consequence of the previous inequality and (4.16). The proof of the claim is complete. ◁

We are now ready to show the validity of (4.13). First we show that by selecting a sufficiently small $\delta>0$, we can guarantee that there is no $x \in K$ with $\beta \Omega(x) \geq 1$, $0<V(x)\leq \delta$ and $|\nabla V(x) f(x)| \leq \delta$. Indeed, by virtue of (4.9) and (4.14), such a $x \in K$ should satisfy the inequalities $\beta \tilde{\gamma}(x) p(|\nabla V(x)|) \geq 1$ and $\sigma(x)\psi(x)|\nabla V(x)|^2 \leq \delta$, which give the inequality $\sigma(x)\psi(x)\left(p^{-1}\left(\dfrac{1}{\beta\tilde{\gamma}(x)}\right)\right)^2 \leq \delta$.

Setting $\delta := \dfrac{1}{2}\min\left\{\sigma(x)\psi(x)\left(p^{-1}\left(\dfrac{1}{\beta\tilde{\gamma}(x)}\right)\right)^2 : x\in K\right\}$ (well-defined and positive since $K$ is compact and since $\tilde{\gamma},\psi,\sigma \in C^0(\Re^n;(0,+\infty))$, $p^{-1} \in K_\infty$), we can guarantee that the inequality $\sigma(x)\psi(x)\left(p^{-1}\left(\dfrac{1}{\beta\tilde{\gamma}(x)}\right)\right)^2 \leq \delta$ does not hold. Consequently, there is no $x \in K$ with $\beta\Omega(x)\geq 1$, $0<V(x)\leq \delta$ and $|\nabla V(x)f(x)|\leq \delta$.

Thus, we are left with the task of showing that for every compact set $K \subseteq \Re^n$ with $K \cap S \neq \varnothing$ there exist a constant $M>0$ with the following property
$$\nabla\theta(x)f(x) \leq \dfrac{M}{\sqrt{V(x)}}|\nabla V(x)f(x)|,$$
for all $x\in K$ with $\beta\Omega(x)<1$, $0<V(x)\leq \delta$ and $|\nabla V(x)f(x)|\leq \delta$ \hfill (4.17)

where $\delta := \dfrac{1}{2}\min\left\{\sigma(x)\psi(x)\left(p^{-1}\left(\dfrac{1}{\beta\tilde{\gamma}(x)}\right)\right)^2 : x\in K\right\}$. Taking into account (4.8) and (4.10), it suffices to show that
$$a(x) \leq \dfrac{M}{\sqrt{V(x)}}\psi(x)|\nabla V(x)|^2,$$
for all $x\in K$ with $\beta\Omega(x)<1$, $0<V(x)\leq \delta$ and $|\nabla V(x)f(x)|\leq \delta$ \hfill (4.18)

Taking into account (4.2), definition (4.5) and the fact that $-|\nabla V(x)|^2|\nabla\theta(x)|^2 + |\nabla V(x)(\nabla\theta(x))'|^2 \leq 0$, we conclude from (4.18) that it suffices to show that
$$-\nabla V(x)(\nabla\theta(x))' \leq \dfrac{M}{\sqrt{\gamma(x)}}|\nabla V(x)|,$$
for all $x \in K$ with $\beta\Omega(x)<1$, $0<V(x)\leq \delta$ and $|\nabla V(x)f(x)|\leq \delta$ \hfill (4.19)

Inequality (4.19) holds by virtue of the Cauchy-Schwarz inequality, if $|\nabla\theta(x)|\sqrt{\gamma(x)}\leq M$. Therefore, the selection $M := 1 + \max\left\{\sqrt{\gamma(x)}|\nabla\theta(x)| : x\in K, \beta\Omega(x)\leq 1\right\}$ is adequate for our purposes. The proof is complete. ◁

When the vector field $F(x)$ can be defined on $\Re^n$ then a simpler formula than the one given (4.6), (4.7) can be used. This is shown in the following result. Its proof is exactly the same with the proof of Theorem 4.1 and is omitted.

**Theorem 4.2:** *Let $V \in C^1(\Re^n; \Re_+)$, $\theta \in C^1(\Re^n; \Re)$ be given functions with locally Lipschitz partial derivatives and suppose that (4.1) holds. Suppose that for every $z \in \Re^n$ and for every $y \in \{x \in \Re^n : V(x) \leq V(z)\}$ the set $\{x \in \Re^n : V(x) \leq V(z), \theta(x) \leq \theta(y)\}$ is compact. Let $\Omega \in C^0(\Re^n; \Re_+)$ be a locally Lipschitz function with $\Omega(x) > 0$ for $x \notin S$, $\Omega(x) = 0$ for all $x \in S$. Suppose that there exist a function $\gamma \in C^0(\tilde{S}; (0, +\infty))$, where $\tilde{S} = \{x \in \Re^n : \Omega(x) < 1\}$, and a locally Lipschitz vector field $F : \Re^n \to \Re^n$ such that (4.2) holds and*

$$\nabla V(x) F(x) \leq 0, \ \nabla \theta(x) F(x) \leq 0, \text{ for all } x \in \Re^n \quad (4.20)$$

*Finally, suppose that there exist locally Lipschitz functions $\psi_i : \Re^n \to (0, +\infty)$ ($i = 1, 2$) for which the following property holds:*

(**) *For every $z \in \Re^n$ there exists $y_z \in \{x \in \Re^n : V(x) \leq V(z)\}$ with $\theta(y_z) \geq \theta(z)$ such that*

$$\{x \in \Re^n : a(x) > 0, V(x) \leq V(z)\} \subseteq \{x \in \Re^n : \theta(x) \leq \theta(y_z)\} \quad (4.21)$$

*where*

$$a(x) := \psi_1(x) \nabla \theta(x) F(x) - \psi_2(x) \nabla V(x) (\nabla \theta(x))' - |\nabla V(x)|^2 |\nabla \theta(x)|^2 + \left|\nabla V(x)(\nabla \theta(x))'\right|^2 \quad (4.22)$$

*Let $\sigma : \Re^n \to (0, +\infty)$ be a locally Lipschitz function and define the locally Lipschitz vector field:*

$$f(x) := \sigma(x)\left(\psi_1(x) F(x) - \psi_2(x)(\nabla V(x))' - \left(|\nabla V(x)|^2 I_n - (\nabla V(x))' \nabla V(x)\right)(\nabla \theta(x))'\right), \text{ for all } x \in \Re^n \quad (4.23)$$

*Then for every $x_0 \in \Re^n$ the unique solution $x(t)$ of the initial value problem (3.4) with $x(0) = x_0$ is defined for all $t \geq 0$ and is bounded. The set $S$, defined in (4.1), is a positively invariant set for the dynamical system (3.4). If we further denote by $\omega(x_0)$ the positive limit set of $\{x(t) : t \geq 0\}$, then for every $x_0 \in \Re^n$, $\omega(x_0)$ is a non-empty, compact, connected, invariant set which satisfies $\omega(x_0) \subseteq \{x \in S : \nabla \theta(x) F(x) = 0\}$. Moreover, every point $x^* \in S$, which satisfies $F(x^*) = 0$, $\theta(x^*) < \theta(x)$ for all $x \in S \setminus \{x^*\}$ with $|x - x^*| < \tilde{\delta}$ and $\{x \in S : \nabla \theta(x) F(x) = 0, |x - x^*| < \tilde{\delta}\} = \{x^*\}$, for an appropriate constant $\tilde{\delta} > 0$, is a locally asymptotically stable equilibrium point of the dynamical system (3.4).*

It should be noted that assumption (**) is less demanding than assumption (*) because the function $a(x)$ defined by (4.22) includes the non-positive term $\nabla \theta(x) F(x)$ (compare with definition (4.5)).

# 5. Construction of a Dynamical NLP Solver

We return to the construction of a feedback law for the control system $\dot{x} = u$ that solves the NLP given by (2.1) and (2.2). As described in the Introduction, the main idea is to use two functions and Theorem 4.2. The first function is a penalty term that penalizes the distance from the feasible set. Here, we will use the penalty function

$$V(x) := \frac{1}{2}|h(x)|^2 + \frac{1}{2}|(g(x))^+|^2 \tag{5.1}$$

Notice that $V \in C^1(\Re^n; \Re_+)$ is a function with locally Lipschitz partial derivatives, since we have $\nabla V(x) = h'(x)A(x) + ((g(x))^+)' B(x)$ for all $x \in \Re^n$, where $A(x), B(x), h(x), g(x)$ are defined in (2.3). However, all what follows can be applied (with appropriate modifications) to functions of the form

$$V(x) := W(h(x)) + \frac{1}{2}\sum_{j=1}^{k} c_j \left(\max(0, g_j(x))\right)^{2p_j}$$

where $W \in C^2(\Re^m; \Re_+)$ is a positive definite, proper function, $c_j > 0$ ($j=1,\ldots,k$) are real constants and $p_j \geq 1$ ($j=1,\ldots,k$) are integers. The second function is the objective function $\theta(x)$.

In order to be able to define an appropriate feedback $u = f(x)$ that guarantees all assumptions of Theorem 3.3, we need the following assumptions.

**(A1)** *For every $z \in \Re^n$ and for every $y \in \{x \in \Re^n : V(x) \leq V(z)\}$ the set $\{x \in \Re^n : V(x) \leq V(z), \theta(x) \leq \theta(y)\}$ is compact.*

Assumption (A1) is a more demanding assumption than (H1).

**(A2)** *For all $x \in S$ the row vectors $\nabla h_i(x)$ ($i = 1,\ldots,m$) and $\nabla g_j(x)$ for all $j = 1,\ldots,k$ for which $g_j(x) = 0$ (active constraints) are linearly independent.*

Assumption (A2) is the linear independence constraint qualification condition. The linear independence constraint qualification, which is assumed in this work is a restrictive assumption: it is more restrictive than the Mangasarian-Fromovitz constraint qualification (assumption (H2)) or the constant rank constraint qualification, which are all more restrictive than the Guignard constraint qualification. However, the linear independence constraint qualification has the advantage of being easily checkable and of being true in many interesting cases (the recent work [25] showed that this assumption holds generically) and it is a vital ingredient for many numerical methods (successive quadratic programming; see for instance [10,26]).

**(A3)** *The following implication holds:*

$$A'(x)h(x) + B'(x)(g(x))^+ = 0 \Rightarrow h(x) = 0 \text{ and } (g(x))^+ = 0 \tag{5.2}$$

*where $A(x), B(x), h(x), g(x)$ are defined in (2.3).*

Assumption (A3) guarantees that there are no critical points of the penalty function out of the feasible set.

Notice that the fact that the symmetric matrix $(A(x)A'(x))$ is positive semidefinite implies $\det(A(x)A'(x)) \geq 0$. Consequently, the condition that the row vectors $\nabla h_i(x)$ ($i=1,...,m$) are linearly independent (being equivalent to $\det(A(x)A'(x)) \neq 0$) is equivalent to the condition $\det(A(x)A'(x)) > 0$.

We next define the symmetric matrix:

$$H(x) = \det(A(x)A'(x))I_n - A'(x)adj(A(x)A'(x))A(x), \text{ for all } x \in \Re^n \qquad (5.3)$$

where the matrix $A(x)$ is defined in (2.3). The following facts are direct consequences of definition (5.3):

**Fact 1:** $H'(x)H(x) = H^2(x) = \det(A(x)A'(x))H(x)$, $A(x)H(x) = 0$ and $H(x)A'(x) = 0$.

**Fact 2:** $\det(A(x)A'(x))\xi'H(x)\xi = |H(x)\xi|^2$, for all $\xi \in \Re^n$

**Fact 3:** For every $\xi \in \Re^n$ and $x \in \Re^n$ with $\det(A(x)A'(x)) > 0$ there exists $\lambda \in \Re^m$ such that

$$\xi = \frac{1}{\det(A(x)A'(x))} H(x)\xi - A'(x)\lambda.$$

We next define the symmetric matrix:

$$Q(x) := \det(A(x)A'(x))B(x)H(x)B'(x) - (\det(A(x)A'(x)))^2 diag((g(x))^-), \text{ for all } x \in \Re^n \qquad (5.4)$$

where $B(x), g(x)$ are defined in (2.3) and $H(x)$ is defined in (5.3). Again the matrix $Q(x) \in \Re^{k \times k}$ is positive semidefinite, since by virtue of Fact 2, the following equality holds for all $\xi = (\xi_1,...,\xi_k)' \in \Re^k$:

$$\xi'Q(x)\xi = |H(x)B'(x)\xi|^2 + (\det(A(x)A'(x)))^2 \sum_{j=1}^{k} |\min(0, g_j(x))|\xi_j^2 \qquad (5.5)$$

Therefore, we get:

$$\det(Q(x)) \geq 0, \text{ for all } x \in \Re^n \qquad (5.6)$$

The following lemma provides necessary and sufficient conditions for the matrix $Q(x) \in \Re^{k \times k}$ to be positive definite. Its proof is provided in the Appendix.

**Lemma 5.1:** *The following statements are equivalent:*
  (a) *The row vectors $\nabla h_i(x)$ ($i=1,...,m$) and $\nabla g_j(x)$ for all $j=1,...,k$ for which $g_j(x) \geq 0$ are linearly independent.*
  (b) *The matrix $Q(x) \in \Re^{k \times k}$ defined by (5.4) is positive definite.*
  (c) $\det(Q(x)) > 0$
  (d) $\det(Q(x)) \neq 0$

Assumption (A2) allows us to construct a vector field $F(x)$ for all $x \in \Re^n$, which satisfies $\nabla \theta(x)F(x) \leq 0$ and $\nabla V(x)F(x) \leq 0$ for all $x \in \Re^n$. This is achieved by the following lemma. Its proof is provided in the Appendix.

**Lemma 5.2:** *Suppose that assumption (A2) holds. Let $R(x) \in \Re^{m \times k}$ be the matrix defined by*

$$R(x) := H(x)B'(x)adj(Q(x)) \quad (5.7)$$

*where $H(x)$ is defined in (5.3), $Q(x) \in \Re^{k \times k}$ is defined in (5.4) and $B(x) \in \Re^{k \times n}$ is defined in (2.3). Let $F(x) \in \Re^n$ be the vector defined by*

$$F(x) := (\det(A(x)A'(x)))^4 R(x) diag((g(x))^-) R'(x)(\nabla \theta(x))' - (\det(A(x)A'(x)))^4 R(x)\left[R'(x)(\nabla \theta(x))'\right]^+$$
$$- \det(A(x)A'(x))(\det(Q(x))I_n - \det(A(x)A'(x))R(x)B(x))H(x)(\det(Q(x))I_n - \det(A(x)A'(x))B'(x)R'(x))(\nabla \theta(x))' \quad (5.8)$$

*Then the following inequalities hold:*

$$\nabla \theta(x) F(x) = (\det(A(x)A'(x)))^4 \nabla \theta(x) R(x) diag((g(x))^-) R'(x)(\nabla \theta(x))' - (\det(A(x)A'(x)))^4 \left|\left[R'(x)(\nabla \theta(x))'\right]^+\right|^2$$
$$- \left|H(x)(\det(Q(x))I_n - \det(A(x)A'(x))B'(x)R'(x))(\nabla \theta(x))'\right|^2 \leq 0 \quad (5.9)$$

$$\left(h'(x)A(x) + ((g(x))^+)' B(x)\right) F(x) = -(\det(A(x)A'(x)))^3 \det(Q(x))((g(x))^+)' \left[R'(x)(\nabla \theta(x))'\right]^+ \leq 0 \quad (5.10)$$

*Moreover, the following implications hold:*

$$\left.\begin{array}{c} x \in S \\ \nabla \theta(x) F(x) = 0 \end{array}\right\} \Leftrightarrow \left.\begin{array}{c} x \in S \\ F(x) = 0 \end{array}\right\} \Leftrightarrow x \in \Phi \quad (5.11)$$

*where $\Phi \subseteq S$ is the set of Karush-Kuhn-Tucker (KKT) points for the problem defined by (2.1) and (2.2).*

For our purposes, we also need a locally Lipschitz function $\Omega \in C^0(\Re^n; \Re_+)$ with $\Omega(x) > 0$ for $x \notin S$, $\Omega(x) = 0$ for all $x \in S$ and such that the following implication hold

$$\Omega(x) < 1 \Rightarrow \det(Q(x)) > 0 \quad (5.12)$$

where $Q(x) \in \Re^{k \times k}$ is defined by (5.4). Such a function can be found easily. For example, the function

$$\Omega(x) := \frac{(1 + c_1(x))V(x)}{c_2(x)\det(Q(x)) + V(x)} \quad (5.13)$$

where $c_i \in C^0(\Re^n; (0, +\infty))$ ($i = 1, 2$) are arbitrary locally Lipschitz functions, satisfies implication (5.12) as well as the requirements $\Omega(x) > 0$ for $x \notin S$, $\Omega(x) = 0$ for all $x \in S$. Moreover, by virtue of assumption (A2) and Lemma 5.1, $\Omega$ as given by (5.13) is defined on $\Re^n$ and is a locally Lipschitz function.

We are now in a position to give our result for the dynamical NLP solver. The result is based on Theorem 4.2.

**Theorem 5.3:** *Suppose that assumptions (A1), (A2), (A3) hold for the NLP problem defined by (2.1) and (2.2) as well as the following assumption:*

**(A4)** *There exist locally Lipschitz functions $\psi_i \in C^0(\Re^n;(0,+\infty))$ ($i=1,2$) such that the following property holds: for every $z \in \Re^n$ there exists $y_z \in \{x \in \Re^n : V(x) \leq V(z)\}$ with $\theta(y_z) \geq \theta(z)$ such that*

$$\{x \in \Re^n : a(x) > 0, V(x) \leq V(z)\} \subseteq \{x \in \Re^n : \theta(x) \leq \theta(y_z)\} \quad (5.14)$$

*where*

$$a(x) := \psi_1(x)\nabla\theta(x)F(x) - |\nabla V(x)|^2|\nabla\theta(x)|^2 - \psi_2(x)\nabla V(x)(\nabla\theta(x))' + |\nabla V(x)(\nabla\theta(x))'|^2 \quad (5.15)$$

*and $V \in C^1(\Re^n; \Re_+)$ is the function defined by (5.1), $F : \Re^n \to \Re^n$ is the locally Lipschitz vector field defined by (5.7), (5.8). Let $\sigma : \Re^n \to (0,+\infty)$ be an arbitrary locally Lipschitz function. Define the locally Lipschitz vector field:*

$$f(x) := \sigma(x)\psi_1(x)F(x) - \sigma(x)\left(\psi_2(x)(\nabla V(x))' + \left(|\nabla V(x)|^2 I_n - (\nabla V(x))'\nabla V(x)\right)(\nabla\theta(x))'\right), \text{ for } x \in \Re^n \quad (5.16)$$

*Let $\Phi \subseteq S$ be the set of KKT points for the problem defined by (2.1) and (2.2). Then the following properties hold for the dynamical system (3.4):*

i) *For every $x_0 \in \Re^n$ the unique solution $x(t)$ of the initial value problem (3.4) with $x(0) = x_0$ is defined for all $t \geq 0$ and is bounded. Moreover, $\omega(x_0)$ is a non-empty, compact, connected, invariant set which satisfies $\omega(x_0) \subseteq \Phi$.*
ii) *Every KKT point of the NLP problem described by (2.1) and (2.2) is an equilibrium point of the dynamical system (3.4) and every equilibrium point of the dynamical system (3.4) is a KKT point of the NLP problem described by (2.1) and (2.2).*
iii) *Every isolated KKT point, which is a strict local minimum of the NLP problem described by (2.1) and (2.2) is a locally asymptotically stable equilibrium point of the dynamical system (3.4).*
iv) *The feasible set $S$, defined in (2.2), is a positively invariant set for the dynamical system (3.4).*

**Proof:** We use Theorem 4.2 for the function $V$ defined by (5.1). The conclusions of the theorem are direct consequences of Theorem 4.2 and Lemma 5.2.

Notice that Assumption (A3) guarantees that (4.1) holds. Next, we show that there exists $\gamma \in C^0(\tilde{S};(0,+\infty))$ such that $\gamma(x)|\nabla V(x)|^2 \geq V(x)$ holds for all $x \in \Re^n$ with $\Omega(x) < 1$, where $\tilde{S} = \{x \in \Re^n : \Omega(x) < 1\}$ and $\Omega \in C^0(\Re^n; \Re_+)$ being an arbitrary locally Lipschitz function with $\Omega(x) > 0$ for $x \notin S$, $\Omega(x) = 0$ for all $x \in S$, satisfying implication (5.12) (e.g., the function defined in (5.13)).

Since the matrices $(A(x)A'(x)), Q(x)$ are positive definite (see Lemma 5.1 and Assumption (A2)) and continuous on $\tilde{S} = \{x \in \Re^n : \Omega(x) < 1\}$, there exist continuous functions $K_i : \tilde{S} \to (0,+\infty)$ ($i=1,2$) such that

$$\xi'(A(x)A'(x))\xi \geq K_1(x)|\xi|^2, \text{ for all } \xi \in \Re^m, x \in \tilde{S} \quad (5.17)$$

$$\frac{1}{(\det(A(x)A'(x)))^3}\xi'Q(x)\xi \geq K_2(x)|\xi|^2, \text{ for all } \xi \in \Re^k, x \in \tilde{S} \quad (5.18)$$

Notice that definition (5.1) implies that $(\nabla V(x))' = A'(x)h(x) + B'(x)(g(x))^+$, which combined with definition (5.3), can also be written in the following form for all $x \in \tilde{S}$:

$$(\nabla V(x))' = A'(x)\left(h(x) + (A(x)A'(x))^{-1} A(x)B'(x)(g(x))^+\right) + \frac{1}{\det(A(x)A'(x))} H(x)B'(x)(g(x))^+ \quad (5.19)$$

Using (5.19), Facts 1 and 2, we obtain for all $x \in \tilde{S}$:

$$|\nabla V(x)|^2 = \frac{1}{(\det(A(x)A'(x)))^2}\left((g(x))^+\right)' B(x)H(x)B'(x)(g(x))^+ +$$
$$\left(h(x) + (A(x)A'(x))^{-1} A(x)B'(x)(g(x))^+\right)'(A(x)A'(x))\left(h(x) + (A(x)A'(x))^{-1} A(x)B'(x)(g(x))^+\right) \quad (5.20)$$

Using (5.20) and the fact that $\left((g(x))^+\right)' diag\left((g(x))^-\right)(g(x))^+ = 0$ in conjunction with definition (5.4), we obtain for all $x \in \tilde{S}$:

$$|\nabla V(x)|^2 = \frac{1}{(\det(A(x)A'(x)))^3}\left((g(x))^+\right)' Q(x)(g(x))^+ +$$
$$\left(h(x) + (A(x)A'(x))^{-1} A(x)B'(x)(g(x))^+\right)'(A(x)A'(x))\left(h(x) + (A(x)A'(x))^{-1} A(x)B'(x)(g(x))^+\right)$$

or

$$|\nabla V(x)|^2 = \left((g(x))^+\right)'\left(\frac{1}{(\det(A(x)A'(x)))^3} Q(x) + B(x)A'(x)(A(x)A'(x))^{-1} A(x)B'(x)\right)(g(x))^+ \quad (5.21)$$
$$+ h'(x)(A(x)A'(x))h(x) + 2h'(x)A(x)B'(x)(g(x))^+$$

Using (5.1), (5.17), (5.18), (5.21), the inequality $-z'Gz - y'G^{-1}y \leq 2z'y$ (which holds for every $y, z \in \Re^m$ and for every positive definite matrix $G \in \Re^{m \times m}$) with $z = h(x)$, $G = \varepsilon(x)(A(x)A'(x))$, $y = A(x)B'(x)(g(x))^+$, where $\varepsilon(x) := \frac{q(x)+1}{K_2(x)+q(x)+1}$ and $q: \tilde{S} \to (0,+\infty)$ is any continuous function that satisfies $\left|B(x)A'(x)(A(x)A'(x))^{-1} A(x)B'(x)\right| \leq q(x)$ for all $x \in \tilde{S}$, we obtain for all $x \in \tilde{S}$:

$$|\nabla V(x)|^2 \geq \left((g(x))^+\right)'\left(\frac{1}{(\det(A(x)A'(x)))^3} Q(x) + (1-\varepsilon^{-1}(x))B(x)A'(x)(A(x)A'(x))^{-1} A(x)B'(x)\right)(g(x))^+$$
$$+ (1-\varepsilon(x))h'(x)(A(x)A'(x))h(x)$$
$$\geq \frac{1}{(\det(A(x)A'(x)))^3}\left((g(x))^+\right)' Q(x)(g(x))^+ + \frac{K_2(x)}{K_2(x)+q(x)+1} h'(x)(A(x)A'(x))h(x)$$
$$- \frac{K_2(x)}{q(x)+1}\left((g(x))^+\right)' B(x)A'(x)(A(x)A'(x))^{-1} A(x)B'(x)(g(x))^+$$
$$\geq K_2(x)\left|(g(x))^+\right|^2 + \frac{K_2(x)K_1(x)}{K_2(x)+q(x)+1}|h(x)|^2 - \frac{K_2(x)}{q(x)+1}\left|B(x)A'(x)(A(x)A'(x))^{-1} A(x)B'(x)\right|\left|(g(x))^+\right|^2$$
$$\geq K_2(x)\left|(g(x))^+\right|^2 + \frac{K_2(x)K_1(x)}{K_2(x)+q(x)+1}|h(x)|^2 - \frac{K_2(x)}{q(x)+1} q(x)\left|(g(x))^+\right|^2$$
$$= K_2(x)\frac{1}{q(x)+1}\left|(g(x))^+\right|^2 + \frac{K_2(x)K_1(x)}{K_2(x)+q(x)+1}|h(x)|^2$$
$$\geq \frac{K_2(x)\min(1, K_1(x))}{K_2(x)+q(x)+1}\left(\left|(g(x))^+\right|^2 + |h(x)|^2\right) = 2\frac{K_2(x)\min(1, K_1(x))}{K_2(x)+q(x)+1} V(x)$$

The above inequality shows that there exists $\gamma \in C^0(\tilde{S};(0,+\infty))$ such that $\gamma(x)|\nabla V(x)|^2 \geq V(x)$ holds for all $x \in \Re^n$ with $\Omega(x) < 1$ (e.g., $\gamma(x) = \dfrac{K_2(x) + q(x) + 1}{2K_2(x)\min(1, K_1(x))}$).

All rest assumptions of Theorem 4.2 are direct consequences of Assumptions (A1), (A2), (A3), (A4) and Lemma 5.2. The proof is complete. ◁

To understand that the proposed NLP solver is an extension of "steepest descent" NLP solvers for unconstrained problems, we can consider the unconstrained NLP problem (2.1) with $S := \Re^n$. In order to apply Theorem 5.3, we can perform the following steps:

i) We can add the scalar inequality constraint $g(x) \equiv -1$ ($k = 1$).

ii) We can add one more state variable $x_{n+1}$ and the scalar equality constraint $x_{n+1} = 0$ ($m = 1$).

Computing the vector field $F(x)$ defined by (5.8), it becomes clear that all assumptions (A1)-(A4) hold with arbitrary locally Lipschitz functions $\psi_i \in C^0(\Re^n;(0,+\infty))$ ($i = 1, 2$), provided that the following assumption holds.

**(A5)** *For every $y \in \Re^n$, the set $\{x \in \Re^n : \theta(x) \leq \theta(y)\}$ is compact.*

Therefore, computing the vector field $f(x)$ defined by (5.16), we are in a position to obtain the following corollary.

**Corollary 5.4:** *Suppose that assumption (A5) holds for the NLP problem defined by (2.1) with $S := \Re^n$. Let $\sigma : \Re^n \to (0,+\infty)$ be an arbitrary locally Lipschitz function. Define the locally Lipschitz vector field:*

$$f(x) := -\sigma(x)(\nabla \theta(x))', \text{ for } x \in \Re^n \qquad (5.22)$$

*Let $\Phi = \{x \in \Re^n : \nabla \theta(x) = 0\}$ (the set of critical points for the problem defined by (2.1) with $S := \Re^n$). Then the following properties hold for the dynamical system (3.4):*
  i) *For every $x_0 \in \Re^n$ the unique solution $x(t)$ of the initial value problem (3.4) with $x(0) = x_0$ is defined for all $t \geq 0$ and is bounded. Moreover, $\omega(x_0)$ is a non-empty, compact, connected, invariant set which satisfies $\omega(x_0) \subseteq \Phi$.*
  ii) *Every critical point of the NLP problem described by (2.1) with $S := \Re^n$ is an equilibrium point of the dynamical system (3.4) and every equilibrium point of the dynamical system (3.4) is a critical point of the NLP problem described by (2.1) with $S := \Re^n$.*
  iii) *Every isolated critical point, which is a strict local minimum of the NLP problem described by (2.1) with $S := \Re^n$ is a locally asymptotically stable equilibrium point of the dynamical system (3.4).*

The conclusions of Corollary 5.4 are almost trivial. Corollary 5.4 is not stated here for its usefulness but for another reason: Corollary 5.4 shows that the NLP solver constructed by Theorem 5.3 is a direct extension of "steepest descent" NLP solvers for unconstrained problems.

# 6. Special Cases

In this section we provide simpler formulas for certain special cases.

1st Case: No equality constraints.

In this case $S := \left\{ x \in \Re^n : \max_{j=1,\ldots k} (g_j(x)) \leq 0 \right\}$. For this case, we add one more state variable $x_{n+1}$ and the equality constraint $x_{n+1} = 0$. Performing all calculations of Theorem 5.3, we are in a position to show that the dynamical NLP solver works under the following assumptions:

**(A1')** *For every $z \in \Re^n$ and for every $y \in \left\{ x \in \Re^n : |(g(x))^+| \leq |(g(z))^+| \right\}$, where $g(x)$ is defined by (2.3), the set $\left\{ x \in \Re^n : |(g(x))^+| \leq |(g(z))^+|, \theta(x) \leq \theta(y) \right\}$ is compact.*

**(A2')** *For all $x \in S$ the row vectors $\nabla g_j(x)$ for all $j = 1,\ldots,k$ for which $g_j(x) = 0$ are linearly independent.*

**(A3')** *The following implication holds: $B'(x)(g(x))^+ = 0 \Rightarrow (g(x))^+ = 0$, where $B(x), g(x)$ are defined by (2.3).*

**(A4')** *There exist locally Lipschitz functions $\psi_i \in C^0(\Re^n; (0,+\infty))$ ($i = 1,2$) such that the following property holds: for every $z \in \Re^n$ there exists $y_z \in \left\{ x \in \Re^n : |(g(x))^+| \leq |(g(z))^+| \right\}$ with $\theta(y_z) \geq \theta(z)$ such that*

$$\left\{ x \in \Re^n : a(x) > 0, |(g(x))^+| \leq |(g(z))^+| \right\} \subseteq \left\{ x \in \Re^n : \theta(x) \leq \theta(y_z) \right\}$$

*where*

$$a(x) := \psi_1(x) \nabla \theta(x) F(x) - |B'(x)(g(x))^+|^2 |\nabla \theta(x)|^2 + (\nabla \theta(x) B'(x)(g(x))^+)^2 - \psi_2(x) \nabla \theta(x) B'(x)(g(x))^+$$

$$Q(x) := B(x)B'(x) - diag((g(x))^-), \quad R(x) := B'(x) adj(Q(x)) \qquad (6.1)$$

$$F(x) := R(x) diag((g(x))^-) R'(x)(\nabla \theta(x))' - R(x) \left( R'(x)(\nabla \theta(x))' \right)^+$$
$$- (\det(Q(x))I_n - R(x)B(x))(\det(Q(x))I_n - B'(x)R'(x))(\nabla \theta(x))' \qquad (6.2)$$

In this case the proposed dynamical NLP solver is defined for every locally Lipschitz function $\sigma : \Re^n \to (0,+\infty)$ by the formula:

$$f(x) := \sigma(x)\left( \psi_1(x) F(x) - \psi_2(x) B'(x)(g(x))^+ - \left( |B'(x)(g(x))^+|^2 I_n - B'(x)(g(x))^+ ((g(x))^+)' B(x) \right)(\nabla \theta(x))' \right),$$

$$\text{for } x \in \Re^n \qquad (6.3)$$

where $F(x)$ is defined by (6.1), (6.2) and $B(x), g(x)$ are defined in (2.3). Notice that assumptions (A1'), (A4') hold automatically for arbitrary locally Lipschitz functions $\psi_i \in C^0(\Re^n; (0,+\infty))$ ($i = 1,2$), if the sets $S_c := \left\{ x \in \Re^n : \max_{j=1,\ldots k} (g_j(x)) \leq c \right\}$ are compact for every $c \geq 0$.

2nd Case: No inequality constraints.

In this case $S := \left\{ x \in \Re^n : h_1(x) = \ldots = h_m(x) = 0 \right\}$. For this case, we add the inequality constraint $g(x) \equiv -1$ and we are in a position to show that the dynamical NLP solver works under the following assumptions:

**(A1'')** *For every $z \in \Re^n$ and for every $y \in \left\{ x \in \Re^n : |h(x)| \leq |h(z)| \right\}$, where $h(x)$ is defined by (2.3), the set $\left\{ x \in \Re^n : |h(x)| \leq |h(z)|, \theta(x) \leq \theta(y) \right\}$ is compact.*

**(A2'')** *For all* $x \in S$ *the row vectors* $\nabla h_i(x)$ ( $i = 1,...,m$ ) *are linearly independent, i.e.,* $\det(A(x)A'(x)) > 0$, *where* $A(x)$ *is defined by (2.3).*

**(A3'')** *The following implication holds:* $A'(x)h(x) = 0 \Rightarrow h(x) = 0$.

**(A4'')** *There exist locally Lipschitz functions* $\psi_i \in C^0(\Re^n;(0,+\infty))$ ( $i = 1,2$ ) *such that the following property holds: for every* $z \in \Re^n$ *there exists* $y_z \in \{x \in \Re^n : |h(x)| \le |h(z)|\}$ *with* $\theta(y_z) \ge \theta(z)$ *such that*

$$\{x \in \Re^n : a(x) > 0, |h(x)| \le |h(z)|\} \subseteq \{x \in \Re^n : \theta(x) \le \theta(y_z)\}$$

*where* $H(x)$ *is defined by (5.3) and* $h(x), A(x)$ *are defined by (2.3) and*

$$a(x) := -\psi_1(x)|H(x)(\nabla\theta(x))'|^2 - |A'(x)h(x)|^2|\nabla\theta(x)|^2 + (\nabla\theta(x)A'(x)h(x))^2 - \psi_2(x)\nabla\theta(x)A'(x)h(x) \quad (6.4)$$

In this case the proposed dynamical NLP solver may be defined for an arbitrary locally Lipschitz function $\sigma : \Re^n \to (0,+\infty)$ by the formula:

$$f(x) := -\sigma(x)\left(\psi_1(x)\det(A(x)A'(x))H(x) + |A'(x)h(x)|^2 I_n - A'(x)h(x)h'(x)A(x)\right)(\nabla\theta(x))' - \sigma(x)\psi_2(x)A'(x)h(x),$$
$$\text{for } x \in \Re^n \quad (6.5)$$

where $H(x)$ is defined by (5.3) and $h(x), A(x)$ are defined by (2.3).

## 7. Examples

In order to demonstrate the strength of the obtained results we have used two examples from [35] and one example with a linear equality constraint.

**Example 7.1:** The first example is dealing with the solution of the problem:

$$\min x_1^2 + 2x_2^2 + x_1 x_2 - 6x_1 - 2x_2 - 12x_3$$
$$s.t.$$
$$x_1 + x_2 + x_3 - 2 = 0$$
$$\begin{bmatrix} -x_1 + 2x_2 - 3 \\ -x_1 \\ -x_2 \\ -x_3 \end{bmatrix} \le 0 \quad (7.1)$$

The problem can be turned to a problem with inequality constraints by eliminating $x_3$. We prefer to eliminate $x_3$ because the dynamics of the dynamical NLP solver will be visible from the phase diagram. By eliminating the variable $x_3$, we obtain the following NLP problem:

$$\min \theta(x) = x_1^2 + 2x_2^2 + x_1 x_2 + 6x_1 + 10x_2$$
$$s.t. \quad (7.2)$$
$$g(x) = \begin{bmatrix} -x_1 + 2x_2 - 3 \\ -x_1 \\ -x_2 \\ x_1 + x_2 - 2 \end{bmatrix} \le 0$$

We notice that assumptions (A1'), (A4') hold automatically for arbitrary locally Lipschitz functions $\psi_i \in C^0(\Re^n;(0,+\infty))$ ( $i = 1,2$ ) since the sets $S_c := \left\{x \in \Re^2 : \max_{j=1,...,4}(g_j(x)) \le c\right\}$ are compact for every $c \ge 0$. Moreover, assumptions (A2') and (A3') hold for this problem, as it can be verified by direct calculations. We used formulas (6.1), (6.2), (6.3) with $\sigma(x)\equiv 1$, $\psi_i(x)\equiv 1$ ($i=1,2$) for the construction of the dynamical NLP solver. The phase diagram of the dynamical NLP solver is shown in Figure 1.

The phase diagram in Figure 1 shows global attractivity at $0\in\Re^2$, which was expected from Theorem 5.3 and the fact that for the NLP problem (7.2) we have $\Phi=\{0\in\Re^2\}$. Since $0\in\Re^2$ is a strict local minimum of the NLP problem (7.1), we can conclude that $0\in\Re^2$ is globally asymptotically stable.

The reader may criticize the efficiency of the dynamical NLP solver, since the phase diagram in Fig.1 shows that many trajectories are "sent" to the third quadrant, while the solution is at zero. This happens because some of the trajectories are inevitably attracted for an initial transient period by the (unconstrained) minimizer of the function $\theta(x)=x_1^2+2x_2^2+x_1x_2+6x_1+10x_2$, which is at $(-2,-2)$. During this short initial transient period, a simultaneous decrease of the values of both the objective function $\theta(x)$ and the penalty function $V(x)$ (defined by (5.1) with $h(x)\equiv 0$), is achieved. However, the trajectories are subsequently pushed towards the feasible set. This is shown in Figure 2, where the time evolution of the values of the objective function $\theta(x)$ and the penalty function $50V(x)$ (with $V(x)$ defined by (5.1) with $h(x)\equiv 0$) is plotted for the solution of the dynamical NLP solver (6.1), (6.2), (6.3) with $\sigma(x)\equiv 1$, $\psi_i(x)\equiv 1$ ($i=1,2$) and initial condition $(1.5,-0.5)$. The objective function obtains very quickly negative values and the solution is subsequently pushed smoothly towards the feasible set (which leads to an eventual increase of the value of the objective function). Therefore, there is no "overshoot" in the value of the objective function (the minimizer is approached from below).  ◁

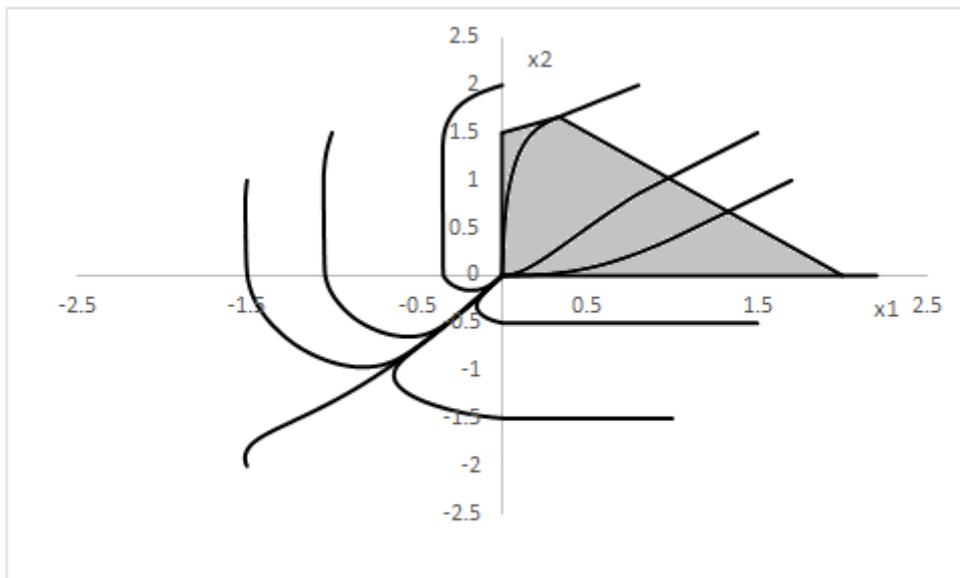

**Fig.1:** The phase diagram of the dynamical NLP solver of Example 7.1. The feasible region is shaded with grey color. Some trajectories are "sent" to the third quadrant, because the trajectories are attracted for an initial transient period by the (unconstrained) minimizer of the objective function $\theta(x)=x_1^2+2x_2^2+x_1x_2+6x_1+10x_2$, which is at $(-2,-2)$. During this short initial transient period, a simultaneous decrease of the values of both the objective function $\theta(x)$ and the penalty function $V(x)$ (defined by (5.1) with $h(x)\equiv 0$), is achieved. However, the trajectories are subsequently pushed towards the feasible set.

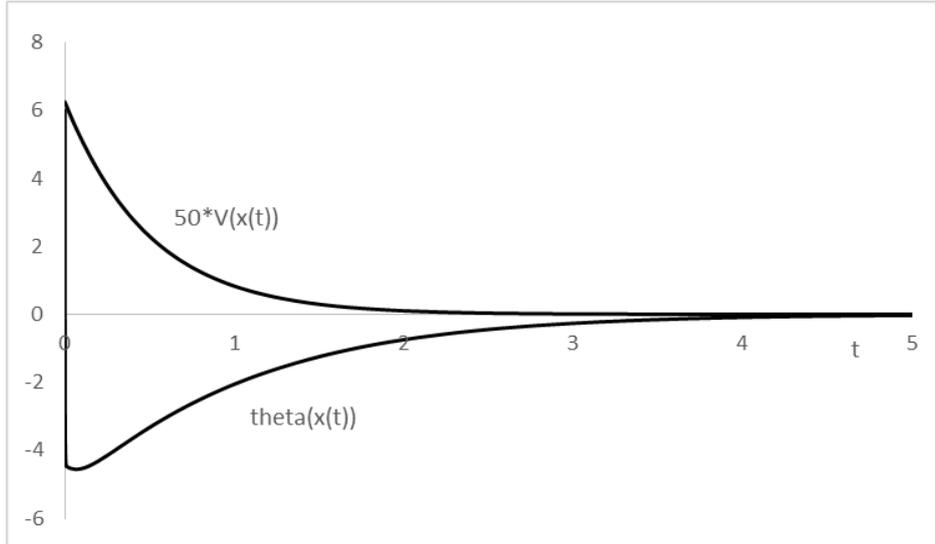

**Fig.2:** The time evolution of the values of the objective function $\theta(x)$ and the penalty function $50V(x)$ ($V(x)$ defined by (5.1) with $h(x) \equiv 0$) for the solution of the dynamical NLP solver (6.1), (6.2), (6.3) with $\sigma(x) \equiv 1$, $\psi_i(x) \equiv 1$ ($i = 1,2$) and initial condition $(1.5, -0.5)$.

**Example 7.2:** Consider the NLP problem

$$\min \theta(x) = x_1^2 + ax_2^2$$
$$\text{s.t.} \qquad (7.3)$$
$$h(x) = x_1 - b = 0$$

with no inequality constraints, where $a > 0$, $b > 0$ are constants. Assumption (A1'') holds for this problem, since the objective function $\theta(x)$ is radially unbounded. Assumptions (A2'') and (A3'') hold trivially (notice that $A(x) = [1 \ 0]$) and it holds that $\det(A(x)A'(x)) > 0$ for all $x \in \Re^2$.

Next we show that assumption (A4'') holds with $\psi_1(x) \equiv 1$ and $\psi_2 \in C^0(\Re^2;(0,+\infty))$ being any positive locally Lipschitz function of only one variable $x_1$, i.e., $\psi_2(x) = \psi_2(x_1)$. We have from (6.4) and (5.3):

$$H(x) = \begin{bmatrix} 0 & 0 \\ 0 & 1 \end{bmatrix}, \ a(x) = -4a^2 x_2^2 \left(1 + (x_1 - b)^2\right) - 2x_1(x_1 - b)\psi_2(x_1) \qquad (7.4)$$

It follows that

$$a(x) > 0 \Leftrightarrow -2a^2 x_2^2 \left(1 + (x_1 - b)^2\right) > x_1(x_1 - b)\psi(x_1) \qquad (7.5)$$

Inequality (7.5) implies that $x_1(x_1 - b) < 0$, or equivalently $x_1 \in (0, b)$. Moreover, inequality (7.5) implies that

$$\theta(x) = x_1^2 + ax_2^2 < x_1^2 - \frac{x_1(x_1 - b)}{2a\left(1 + (x_1 - b)^2\right)}\psi(x_1) \qquad (7.6)$$

For every $x \in \Re^2$ with $|x_1 - b| \leq |z_1 - b|$ that satisfies (7.5) we get from (7.6):

$$\theta(x) \leq C = \max\left\{ x_1^2 - \frac{x_1(x_1 - b)}{2a\left(1 + (x_1 - b)^2\right)}\psi(x_1) : 0 \leq x_1 \leq b \right\} \qquad (7.7)$$

Let $z = (z_1, z_2) \in \Re^2$ be an arbitrary given vector. Inequality (7.7) implies the existence of a vector $y = (y_1, y_2) \in \Re^2$ with $|y_1 - b| \leq |z_1 - b|$ and $\theta(y) \geq \theta(z)$ for which the following implication holds:

$$-2a^2 x_2^2 \left(1 + (x_1 - b)^2\right) > x_1(x_1 - b)\psi(x_1) \quad and \quad |x_1 - b| \leq |z_1 - b| \Rightarrow \theta(x) \leq \theta(y) \qquad (7.8)$$

Such a vector $y = (y_1, y_2) \in \Re^2$ can always be found (e.g. take $y_1 = b$, $y_2 = \sqrt{\dfrac{C - b^2}{a}}$ if $C \geq \theta(z)$ and $y = z$ if $C < \theta(z)$ ). Therefore, assumption (A4'') holds.

For this problem, the dynamical NLP solver (6.5) for $\psi_1(x) \equiv 1$ and arbitrary locally Lipschitz functions $\sigma : \Re^2 \to (0,+\infty)$, $\psi_2 \in C^0(\Re;(0,+\infty))$ is given by:

$$f(x) := -\sigma(x) \begin{bmatrix} \psi_2(x_1)(x_1 - b) \\ 2ax_2\left(1 + (x_1 - b)^2\right) \end{bmatrix}, \text{ for } x \in \Re^2 \tag{7.9}$$

Using the Lyapunov function $W(x) = (x_1 - b)^2 / 2 + x_2^2 / 2$, it can be shown that the equilibrium point $(b,0)' \in \Re^2$ of system $\dot{x} = f(x)$ with $\sigma(x) \equiv c_1 > 0$, $\psi_2(x) \equiv c_2 > 0$ is globally exponentially stable (see [21]). This is a stronger property than the global asymptotic stability property, which was expected from Theorem 5.3, the fact that $(b,0)' \in \Re^2$ is a strict local minimum of the NLP problem (7.3) and the fact that for the NLP problem (7.3) we have $\Phi = \{(b,0)' \in \Re^2\}$. ◁

**Example 7.3:** The third example is the Rosen–Suzuki problem:

$$\min \theta(x) = x_1^2 + x_2^2 + 2x_3^2 + x_4^2 - 5x_1 - 5x_2 - 21x_3 + 7x_4$$

s.t.

$$h(x) = 2x_1^2 + x_2^2 + x_3^2 + 2x_1 - x_2 - x_4 - 5 = 0 \tag{7.10}$$

$$g(x) = \begin{bmatrix} x_1^2 + x_2^2 + x_3^2 + x_4^2 + x_1 - x_2 + x_3 - x_4 - 8 \\ x_1^2 + 2x_2^2 + x_3^2 + 2x_4^2 - x_1 - x_4 - 10 \end{bmatrix} \leq 0$$

It should be noticed that (7.10) is an NLP problem with nonlinear equality and inequality constraints. For this problem we have $\det(A(x)A'(x)) > 0$ for all $x \in \Re^4$. We notice that assumptions (A1), (A4) hold automatically for every locally Lipschitz functions $\psi_i \in C^0(\Re^n;(0,+\infty))$ ($i = 1,2$) since the sets $S_c := \left\{ x \in \Re^4 : |h(x)| \leq c, \max_{j=1,\ldots,4}(g_j(x)) \leq c \right\}$ are compact for every $c \geq 0$. Moreover, assumptions (A2) and (A3) hold for this problem, as it can be verified by direct (but tedious) calculations. We used formula (5.25) with

$$\sigma(x) := \frac{1}{1 + \left|\psi_1(x)F(x) - (\nabla V(x))' - \left(|\nabla V(x)|^2 I_n - (\nabla V(x))'\nabla V(x)\right)(\nabla \theta(x))'\right|}, \quad \psi_1(x) = \frac{1}{(\det(A(x)A'(x)))^{10}}, \quad \psi_2(x) \equiv 1$$

for the construction of the dynamical NLP solver. The solution for various initial points are shown in Figure 3. The solutions were obtained by using the subroutine ODE23T in MATLAB.

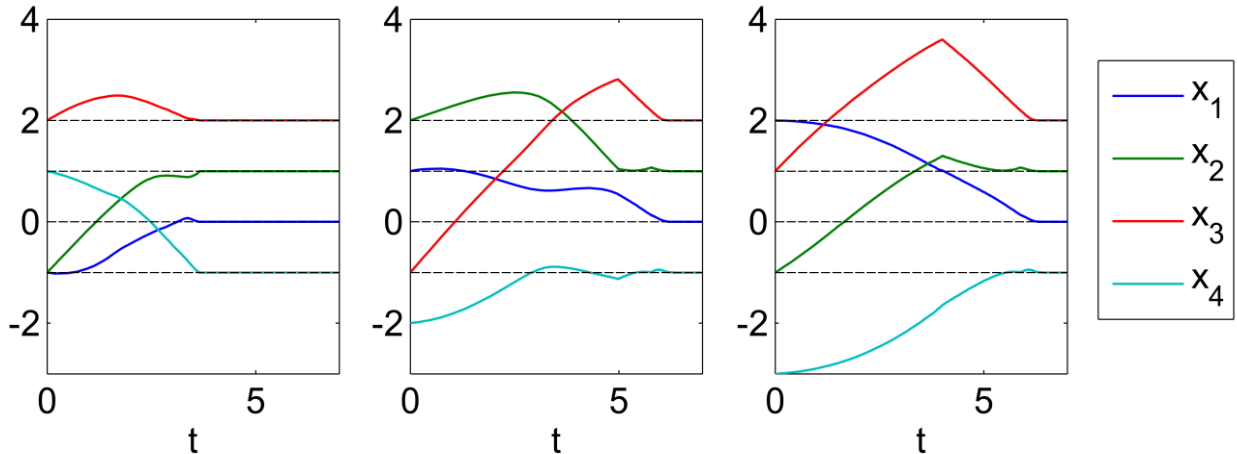

**Fig. 3:** Solution of the dynamical NLP solver of Example 7.3 from various initial points.

In all cases, we observe convergence to $x^* = (0,1,2,-1)' \in \Re^4$. Notice that in Figure 3 one of the initial points is the point $x = (-1,-1,2,1)' \in \Re^4$, which is a special point where the algorithms proposed in [20] could not be used. Even for this initial point, the solution converges rapidly to $x^* = (0,1,2,-1)' \in \Re^4$. ◁

## 8. Concluding Remarks

In this work we have showed that given a nonlinear programming problem, it is possible, under mild assumptions, to construct a family of globally defined dynamical systems, so that: (a) the equilibrium points are the unknown critical points of the problem, (b) for every initial condition, the solution of the corresponding initial value problem converges to the set of critical points, (c) every strict local minimum is locally asymptotically stable, (d) the feasible set is a positively invariant set, and (e) the dynamical system is given explicitly and does not involve the unknown critical points of the problem. No special convexity assumption was employed. The construction of the family of dynamical systems was based on an extension of the Control Lyapunov Function methodology, which employed extensions of the LaSalle's theorem and are of independent interest. Many examples illustrated the obtained results.

At this point the obtained results have nothing to do with extremum seeking (see [12,22]), but may open the way of using different extremum seeking control schemes in the future for constrained problems. Finally, the extension of the obtained results to non-cooperative games for the determination of Nash equilibria may be achieved: this is a future research topic.

**Acknowledgements:** The authors would like to thank Maria Kontorinaki for her help in the simulations of Example 7.3.

# Appendix

**Proof of Theorem 3.1:** Define

$$S := \{ x \in \Re^n : \nabla V(x) f(x) = 0 \} \quad (A.1)$$

Let $x_0 \in \Re^n$ (arbitrary) and consider the unique solution $x(t)$ of the initial value problem (3.4) with $x(0) = x_0$. The solution is defined on $[0, t_{max})$, where $t_{max} \in (0, +\infty]$ is the maximal existence time of the solution. By virtue of (3.1) we get:

$$\frac{d}{dt} V(x(t)) = \nabla V(x(t)) f(x(t)) \leq 0, \text{ for all } t \in [0, t_{max}) \quad (A.2)$$

and consequently, it follows that

$$x(t) \in \{ y \in \Re^n : V(y) \leq V(x_0) \} \text{ for all } t \in [0, t_{max}). \quad (A.3)$$

Let $y_{x_0} \in \Re^n$ be the vector for which $\theta(y_{x_0}) \geq \theta(x_0)$ and the inclusion

$$\{ x \in \Re^n : V(x) \leq V(x_0) \} \subseteq \{ x \in \Re^n : \theta(x) \leq \theta(y_{x_0}) \} \cup \{ x \in \Re^n : \nabla \theta(x) f(x) \leq 0 \} \quad (A.4)$$

holds. Define

$$\theta_{max} := \max \{ \theta(x) : \theta(x) \leq \theta(y_{x_0}), V(x) \leq V(x_0) \}. \quad (A.5)$$

Notice that definition (A.5) is valid since the set

$$K(x_0) := \{ x \in \Re^n : \theta(x) \leq \theta(y_{x_0}), V(x) \leq V(x_0) \} \quad (A.6)$$

is non-empty and compact (the fact that $\theta(y_{x_0}) \geq \theta(x_0)$ implies that $x_0 \in K(x_0)$). We next make the following claim.

**Claim:** $\theta(x(t)) \leq \theta_{max}$, for all $t \in [0, t_{max})$.

**Proof of Claim:** The proof of the claim is made by contradiction. Suppose that there exists $t \in [0, t_{max})$ such that $\theta(x(t)) > \theta_{max}$. Moreover, notice that since $x_0 \in K(x_0)$, the set $\{ s \in [0,t] : x(s) \in K(x_0) \} \neq \emptyset$ is non-empty. Define $T := \sup \{ s \in [0,t] : x(s) \in K(x_0) \}$. Since $K(x_0) \subset \Re^n$ is closed, it follows that $x(T) \in K(x_0)$. Moreover, definitions (A.5), (A.6) imply that $V(x(T)) \leq V(x_0)$, $\theta(x(T)) \leq \theta_{max}$ and that $T < t$. Moreover, since $T := \sup \{ s \in [0,t] : x(s) \in K(x_0) \}$, it follows that $x(s) \notin K(x_0)$ for all $s \in (T,t]$. By virtue of (A.3) and (A.4), we obtain $\frac{d}{ds} \theta(x(s)) = \nabla \theta(x(s)) f(x(s)) \leq 0$, for all $s \in (T,t]$, which implies $\theta(x(t)) \leq \theta(x(T))$. Since $\theta(x(T)) \leq \theta_{max}$, we get $\theta(x(t)) \leq \theta_{max}$; a contradiction. The proof of the Claim is complete. ◁

Let $y \in K(x_0)$ be such that $\theta(y) = \theta_{max}$. Then by virtue of the Claim and (A.3), we obtain $x(t) \in \{x \in \Re^n : V(x) \le V(x_0), \theta(x) \le \theta(y)\}$ for all $t \in [0, t_{max})$. Since the set $\{x \in \Re^n : V(x) \le V(x_0), \theta(x) \le \theta(y)\}$ is compact, it follows that $t_{max} = +\infty$ (because if $t_{max} < +\infty$ then we should have $\lim_{t \to t_{max}^-} (\|x(t)\|) = +\infty$; a contradiction). Therefore, $x(t)$ is defined for all $t \ge 0$ and is bounded.

LaSalle's theorem (Theorem 3.4 on page 115 in [21]), Theorem 3.3.2.8 on page 120 in [28] and (3.1), (A.1) imply that $\omega(x_0)$ is a non-empty, compact, connected, invariant set which satisfies $\omega(x_0) \subseteq S$.

Application of LaSalle's theorem on the invariant, compact set $\omega(x_0) \subseteq S$ and inequality (3.2) guarantee that for every $y \in \omega(x_0)$ the positive limit set $\omega(y)$ is a non-empty, compact, connected, invariant set which satisfies $\omega(y) \subseteq \{x \in \omega(x_0) : \nabla \theta(x) f(x) = 0\}$. This implies that $\omega(y) \subseteq \omega(x_0) \cap \{x \in \Re^n : \nabla \theta(x) f(x) = \nabla V(x) f(x) = 0\} \ne \emptyset$, for every $y \in \omega(x_0)$.
The proof is complete. ◁

**Proof of Theorem 3.3:** Positive invariance of the set $S$ is a direct consequence of inequality (3.1) and definition (3.5). Let $x_0 \in \Re^n$ be given (arbitrary) and consider the unique solution $x(t)$ of the initial value problem (3.4) with $x(0) = x_0$. By virtue of Theorem 3.1, $\omega(x_0) \subseteq S$ is a compact set.

Let $K \subseteq \Re^n$ be a compact set that contains an open neighborhood $N \subseteq \Re^n$ of $\omega(x_0)$ (such a compact set $K \subseteq \Re^n$ exists since $\omega(x_0) \subseteq S$ is a compact set). Let $g \in K_\infty \cap C^1((0,+\infty); \Re_+)$ and $\delta > 0$ be such that

$$\nabla \theta(x) f(x) \le \frac{dg}{ds}(V(x)) |\nabla V(x) f(x)|, \text{ for all } x \in K \text{ with } 0 < V(x) \le \delta \text{ and } |\nabla V(x) f(x)| \le \delta \quad (A.7)$$

Since $\lim_{t \to +\infty} dist(x(t), \omega(x_0)) = 0$ and since $\omega(x_0)$ is a compact set, it follows from (3.1) and (3.5) that there exists $T > 0$ such that $x(t) \in K$, $\nabla V(x(t)) f(x(t)) \ge -\delta$ and $V(x(t)) \le \delta$ for all $t \ge T$.

Define for all $x \in K$:
$$\tilde{\theta}(x) := \theta(x) + g(V(x)) \quad (A.8)$$

Definition (A.8) in conjunction with inequalities (3.1) and (A.7) implies that:

$$\nabla \tilde{\theta}(x) f(x) \le 0, \text{ for all } x \in K \text{ with } 0 < V(x) \le \delta \text{ and } |\nabla V(x) f(x)| \le \delta \quad (A.9)$$

It follows from (3.2), (3.5), (A.9) (and the fact that $S = \{x \in \Re^n : V(x) = 0\}$ is positively invariant for (3.4); a direct consequence of (3.1)) that the mapping $t \to \tilde{\theta}(x(t))$ is non-increasing for $t \ge T$. Theorem 3.1 implies that the solution $x(t)$ is bounded for all $t \ge 0$. Therefore, the mapping $t \to \tilde{\theta}(x(t))$ is bounded from below. It follows that there exists $l \in \Re$ such that $\lim_{t \to +\infty} \tilde{\theta}(x(t)) = l$. Since $g \in K_\infty$ and $\lim_{t \to +\infty} V(x(t)) = 0$, it follows from definition (A.8) that $\lim_{t \to +\infty} \theta(x(t)) = l$. Therefore, we must have $\omega(x_0) \subseteq \{x \in \Re^n : \theta(x) = l\}$. Invariance of $\omega(x_0) \subseteq S$ implies that $\omega(x_0) \subseteq \{x \in S : \nabla \theta(x) f(x) = 0\}$.

Let $x^* \in S$ be an equilibrium point of the dynamical system (3.4), which satisfies $\theta(x^*) < \theta(x)$ for all $x \in S \setminus \{x^*\}$ with $|x - x^*| < \tilde{\delta}$, $f(x^*) = 0$ and $\{x \in S : \nabla \theta(x) f(x) = 0, |x - x^*| < \tilde{\delta}\} = \{x^*\}$, for an appropriate constant $\tilde{\delta} > 0$. Let $K \subseteq \Re^n$ be a compact set with $\{x \in \Re^n : |x - x^*| < \tilde{\delta}\} \subset K$ and let $g \in K_\infty \cap C^1((0,+\infty); \Re_+)$ be a function satisfying (3.6) for certain constant $\delta > 0$. Consider the function:

$$W(x) = \frac{1}{2}\left(\max\left(0, \theta(x) - \theta(x^*)\right)\right)^2 + g(V(x)) \quad (A.10)$$

defined on the open set

$$D := \left\{ x \in \Re^n : |x - x^*| < \tilde{\delta}, V(x) < \delta, |\nabla V(x) f(x)| < \delta, \theta(x) < \theta(x^*) + 1 \right\} \quad (A.11)$$

Notice that $W : D \to \Re_+$ as defined by (A.10) is continuous. The assumptions for the equilibrium point $x^* \in S$ in conjunction with (3.5) imply that $W(x^*) = 0$ and $W(x) > 0$ for all $x \in D \setminus \{x^*\}$. Definition (A.10) in conjunction with inequalities (3.1) and (3.6) implies that:

$$\nabla W(x) f(x) \leq 0, \text{ for all } x \in D \setminus S \quad (A.12)$$

Positive invariance of $S = \{x \in \Re^n : V(x) = 0\}$ and (3.2) in conjunction with (A.12) and definition (A.10) implies that for every solution $x(t)$ of (3.4) defined on some interval $I \subseteq \Re_+$ and satisfying $x(t) \in D$ for $t \in I$, the mapping $I \ni t \to W(x(t))$ is non-increasing. Theorem 3.3.5 on page 36 in [5] implies that $x^* \in S$ is a stable equilibrium point. Stability implies that there exists $c > 0$ such that the solution $x(t)$ of (3.4) with initial condition $x(0) = x_0$, $|x_0 - x^*| < c$ satisfies $x(t) \in D$ for all $t \geq 0$. Therefore, it follows that $\omega(x_0) \subseteq \{x \in S : \nabla \theta(x) f(x) = 0\} \cap D$ for all $x_0 \in \Re^n$ with $|x_0 - x^*| < c$. Definition (A.11) and the fact that $\{x \in S : \nabla \theta(x) f(x) = 0, |x - x^*| < \tilde{\delta}\} = \{x^*\}$ implies that $\omega(x_0) = \{x^*\}$ for all $x_0 \in \Re^n$ with $|x_0 - x^*| < c$. Consequently, $x^* \in S$ is locally asymptotically stable. The proof is complete. ◁

**Proof of Lemma 5.1:** Equivalence of (b), (c) and (d) is a direct consequence of the fact that the matrix $Q(x) \in \Re^{k \times k}$ defined by (5.4) is positive definite (see (5.5)).

We next show implication (a) $\Rightarrow$ (b) by contradiction. Notice that the linear independence of the row vectors $\nabla h_i(x)$ ($i = 1, \ldots, m$) implies that $\det(A(x) A'(x)) > 0$. Suppose that the matrix $Q(x) \in \Re^{k \times k}$ defined by (5.4) is not positive definite. Then there exists a non-zero $\xi = (\xi_1, \ldots, \xi_k)' \in \Re^k$ with $\xi' Q(x) \xi = 0$. Consequently, equality (5.5) shows that we must have $H(x) B'(x) \xi = 0$ and $\xi_j = 0$ for all $j = 1, \ldots, k$ with $g_j(x) < 0$. Fact 3 implies that there exists $\lambda \in \Re^m$ such that $B'(x) \xi = A'(x) \lambda$. The previous equality implies that

$$\sum_{j=1}^k \xi_j \nabla g_j(x) - \sum_{i=1}^m \lambda_i \nabla h_i(x) = 0 \quad (A.13)$$

Since $\xi_j = 0$ for all $j = 1, \ldots, k$ with $g_j(x) < 0$ and since $\xi = (\xi_1, \ldots, \xi_k)' \in \Re^k$ is non-zero, we conclude from (A.13) that the linear independence of the row vectors $\nabla h_i(x)$ ($i = 1, \ldots, m$) and $\nabla g_j(x)$ for all $j = 1, \ldots, k$ for which $g_j(x) \geq 0$ is violated.

Finally, we show implication (b) $\Rightarrow$ (a) by contradiction. Notice that since $Q(x) = 0$ when $\det(A(x) A'(x)) = 0$, it follows that $\det(A(x) A'(x)) > 0$, or equivalently the row vectors $\nabla h_i(x)$ ($i = 1, \ldots, m$) are linearly independent. Suppose that the row vectors $\nabla h_i(x)$ ($i = 1, \ldots, m$) and $\nabla g_j(x)$ for all $j = 1,\ldots,k$ for which $g_j(x) \geq 0$ are linearly dependent. The linear dependence of the row vectors $\nabla h_i(x)$ ($i = 1,\ldots,m$) and $\nabla g_j(x)$ for all $j = 1,\ldots,k$ for which $g_j(x) \geq 0$ implies the existence of vectors $\xi = (\xi_1,\ldots,\xi_k)' \in \Re^k$, $\lambda \in \Re^m$, not both of them being zero, with $\xi_j = 0$ for all $j = 1,\ldots,k$ with $g_j(x) < 0$, such that (A.13) holds. Linear independence of the row vectors $\nabla h_i(x)$ ($i = 1,\ldots,m$) implies that $\xi = (\xi_1,\ldots,\xi_k)' \in \Re^k$ is not zero. Since (A.13) can be written as $B'(x)\xi = A'(x)\lambda$, it follows from Fact 1 that $H(x)B'(x)\xi = 0$. The facts that $H(x)B'(x)\xi = 0$ and $\xi_j = 0$ for all $j = 1,\ldots,k$ with $g_j(x) < 0$ in conjunction with (5.5) implies that $\xi'Q(x)\xi = 0$, which shows that the matrix $Q(x) \in \Re^{k \times k}$ defined by (5.4) is not positive definite.

The proof is complete. ◁

**Proof of Lemma 5.2:** Relations (5.9), (5.10) are direct consequences of definitions (5.4), (5.7), (5.8), Facts 1, 2 and the facts that $\left((g(x))^+\right)' diag\left((g(x))^-\right) = 0$, $Q(x) adj(Q(x)) = \det(Q(x)) I_k$.

The implications $\left.\begin{array}{l} x \in S \\ \nabla \theta(x) F(x) = 0 \end{array}\right\} \Leftrightarrow \left.\begin{array}{l} x \in S \\ F(x) = 0 \end{array}\right\}$ are direct consequences of definition (5.8) and inequality (5.9) and the fact that assumption (A2) guarantees that $\det(A(x)A'(x)) > 0$ when $x \in S$. We next show the implications $\left.\begin{array}{l} x \in S \\ F(x) = 0 \end{array}\right\} \Leftrightarrow x \in \Phi$.

Suppose that $x \in S, F(x) = 0$. It follows that $\nabla \theta(x) F(x) = 0$. Since assumption (A2) guarantees that $\det(A(x)A'(x)) > 0$ when $x \in S$, we obtain from (5.9)

$$diag\left((g(x))^-\right) R'(x)(\nabla \theta(x))' = 0$$
$$\left(R'(x)(\nabla \theta(x))'\right)^+ = 0$$
$$H(x)\left(\det(Q(x)) I_n - \det(A(x)A'(x)) B'(x) R'(x)\right)(\nabla \theta(x))' = 0$$

or equivalently, since $x \in S$ (which implies that $(g(x))^- = g(x)$)

$$diag(g(x)) R'(x)(\nabla \theta(x))' = 0$$
$$R'(x)(\nabla \theta(x))' \leq 0 \tag{A.14}$$
$$H(x)\left(\det(Q(x)) I_n - \det(A(x)A'(x)) B'(x) R'(x)\right)(\nabla \theta(x))' = 0$$

Lemma 5.1 in conjunction with assumption (A2) implies $\det(Q(x)) > 0$. Define $\mu := -\dfrac{\det(A(x)A'(x))}{\det(Q(x))} R'(x)(\nabla \theta(x))'$ and notice that (A.14) implies that $\mu \geq 0$ and

$$diag(g(x)) \mu = 0 \tag{A.15}$$

Equation (A.15) implies that $\mu' g(x) = 0$. Moreover, (A.14) implies $H(x)\left((\nabla \theta(x))' + B'(x)\mu\right) = 0$. Consequently, Fact 3 implies that there exists $\lambda \in \Re^m$ such that $(\nabla \theta(x))' + B'(x)\mu = -A'(x)\lambda$. Therefore, (2.4) holds and thus $x \in \Phi$.

We finally prove the implication $x \in \Phi \Rightarrow F(x) = 0$. Suppose that there exist $\lambda \in \Re^m$ and $\mu \in \Re^k_+$ such that (2.4) holds. Fact 1 implies that $H(x)\left((\nabla\theta(x))' + B'(x)\mu\right) = 0$. Using (5.7), and the fact that $Q(x), H(x)$ are symmetric matrices (and thus $adj(Q(x))$ is a symmetric matrix), we get

$$R'(x)(\nabla\theta(x))' = adj(Q(x))B(x)H(x)(\nabla\theta(x))' = -adj(Q(x))B(x)H(x)B'(x)\mu \qquad (A.16)$$

Since $\mu \geq 0$ and $g(x) \leq 0$ we also get from $\mu'g(x) = 0$ that (A.15) holds. Combining (A.15), (A.16) and using definition (5.4) and the fact that $x \in \Phi$ implies that $x \in S$ (and thus Lemma 5.1 in conjunction with assumption (A2) implies that $\det(Q(x)) > 0$ and $\det(A(x)A'(x)) > 0$) we get $\mu := -\frac{\det(A(x)A'(x))}{\det(Q(x))} R'(x)(\nabla\theta(x))'$. The facts that $\mu \geq 0$ and $\det(Q(x)) > 0$ imply that $R'(x)(\nabla\theta(x))' \leq 0$, or that $\left(R'(x)(\nabla\theta(x))'\right)^+ = 0$. Moreover, (A.15) in conjunction with the fact that $x \in S$ (which implies that $(g(x))^- = g(x)$) gives $diag\left((g(x))^-\right)R'(x)(\nabla\theta(x))' = 0$. Finally, equation $H(x)\left((\nabla\theta(x))' + B'(x)\mu\right) = 0$ implies that $H(x)\left(\det(Q(x))I_n - \det(A(x)A'(x))B'(x)R'(x)\right)(\nabla\theta(x))' = 0$. Equations $\left(R'(x)(\nabla\theta(x))'\right)^+ = 0$, $diag\left((g(x))^-\right)R'(x)(\nabla\theta(x))' = 0$, $H(x)\left(\det(Q(x))I_n - \det(A(x)A'(x))B'(x)R'(x)\right)(\nabla\theta(x))' = 0$ and definition (5.8) imply that $F(x) = 0$.

The proof is complete. ◁